\input amstex
 \documentstyle{amsppt}
\magnification=\magstep1

%\voffset 1.35 truein
%\hoffset 0.845 truein
\catcode`\@=11
\hsize=125mm
\vsize=184mm
\parskip=0pt

%%FONT SIZE=11pt

\font@\titlebf=cmbx10 scaled \magstep1
\font@\elevrm=cmr10 scaled \magstephalf
\font@\elevit=cmti10 scaled \magstephalf
\font@\elevsl=cmsl10 scaled \magstephalf
\font@\elevbf=cmbx10 scaled \magstephalf
\font@\elevi=cmmi10 scaled \magstephalf
\font@\elevsy=cmsy10 scaled \magstephalf
\font@\elevex=cmex10 scaled \magstephalf
\font@\elevbm=cmmib10 scaled \magstephalf
\font@\eightci=cmmi8
\font@\sixi=cmmi6
\font@\eightcsy=cmsy8
\font@\sixsy=cmsy6
\font@\eightcrm=cmr8
\font@\sixrm=cmr6
\font@\eightcbf=cmbx8
\font@\sixbf=cmbx6
\font@\elevsc=cmcsc10 scaled \magstephalf
\font@\ninerm=cmr9

%%EXTRA SYMBOLS

\newfam\msxfam
\newfam\msyfam
%\newfam\cmmibfam
\font@\elevmsx=msam10 scaled \magstephalf
\font@\eightmsx=msam8
\font@\sixmsx=msam6
\font@\nyolmsx=msam8
\font@\elevmsy=msbm10 scaled \magstephalf
\font@\eightmsy=msbm8
\font@\sixmsy=msbm6
\font@\nyolmsy=msbm8
\def\Bbb{\ifmmode\let\next\Bbb@\else
\def\next{\errmessage{Use \string\Bbb\space only in math mode}}\fi\next}
\def\Bbb@#1{{\Bbb@@{#1}}}
\def\Bbb@@#1{\fam\msyfam#1}

\newfam\tffam
\newtoks\elevpoint@
\def\elevpoint{\normalbaselineskip13.5\p@
\spaceskip3.5\p@ plus2\p@ minus1\p@
\abovedisplayskip13.5\p@ plus3\p@ minus9\p@
\belowdisplayskip\abovedisplayskip
\abovedisplayshortskip\z@ plus3\p@
\belowdisplayshortskip7\p@ plus3\p@ minus4\p@
\textonlyfont@\rm\elevrm \textonlyfont@\it\elevit
\textonlyfont@\sl\elevsl \textonlyfont@\bf\elevbf
\textonlyfont@\tf\elevsc %\textfont\cmmibfam\elevbm
\ifsyntax@ \def\big##1{{\hbox{$\left##1\right.$}}}%
\let\Big\big \let\bigg\big \let\Bigg\big
\else
\textfont\z@\elevrm \scriptfont\z@\eightrm
\scriptscriptfont\z@\sixrm
\textfont\@ne\elevi \scriptfont\@ne\eighti
\scriptscriptfont\@ne\sixi
\textfont\tw@\elevsy \scriptfont\tw@\eightsy
\scriptscriptfont\tw@\sixsy
\textfont\thr@@\elevex \scriptfont\thr@@\eightex
\scriptscriptfont\thr@@\eightex
\textfont\itfam\elevit \scriptfont\itfam\eightit
\scriptscriptfont\itfam\eightit
\textfont\bffam\elevbf \scriptfont\bffam\eightbf
\scriptscriptfont\bffam\sixbf
\textfont\msxfam=\elevmsx \scriptfont\msxfam=\eightmsx
\scriptscriptfont\msxfam=\sixmsx
\textfont\msyfam=\elevmsy \scriptfont\msyfam=\eightmsy
\scriptscriptfont\msyfam=\sixmsy
\setbox\strutbox\hbox{\vrule height8.5\p@ depth3.5\p@ width\z@}%
\setbox\strutbox@\hbox{\lower.5\normallineskiplimit\vbox{%
\kern-\normallineskiplimit\copy\strutbox}}%
\setbox\z@\vbox{\hbox{$($}\kern\z@}\bigsize@1.2\ht\z@
\fi
\normalbaselines\rm\dotsspace@1.5mu\ex@.2326ex\jot3\ex@
\the\elevpoint@}
\elevpoint\rm

%%PAGE OUTPUT

\NoPageNumbers
\def\lefthead{} \def\righthead{}
\headline{\ifnum\pageno=1 \hfill\else\ifodd \pageno
{\phantom{\ninerm \number\pageno}\hfill \ninerm\righthead \hfill
\ninerm\number\pageno} \else {\ninerm\number\pageno \hfill
\ninerm\lefthead\hfill\phantom{\ninerm\number\pageno}}\fi\fi}
\def\headLine#1#2{\def\lefthead{#1}\def\righthead{#2}}
\catcode`\@=12

%%DEFINITIONS

\def\bR{\bold R}

\def\bv{\bold v}
\def\bz{\bold z}
\def\b1{\bold 1}
\def\bmu{\hbox{\elevbm\char'26}}

\def\al{\alpha}

\def\de{\delta}
\def\si{\sigma}

\def\Ga{\Gamma}
\def\De{\Delta}

\def\na{\nabla}

\def\Tr{\operatorname{Tr}}
\def\<{\langle}
\def\>{\rangle}
\def\pa{\partial}

\def\wt{\widetilde}

\topmatter
\title%\titlebf
\centerline{A cornucopia of isospectral pairs of metrics}
\centerline{on spheres
with different local geometries}
\endtitle
\author\tf Z.I. Szab\'o\endauthor
\thanks Research partially supported by
NSF grant DMS-0104361 and 
CUNY grant 9-91907.
\endthanks
\headLine{Z. I. SZAB\'O}{CORNUCOPIA OF ISOSPECTRAL PAIRS}
\address
Zolt\'an Imre Szab\'o
\endaddress
\address
\noindent
City University of New York,
Lehman College,
Bronx, NY
\endaddress
\address
R\'enyi Alfr\'ed Institute of Mathematics,
Budapest,
Hungary
\endaddress
\address
Email: {\rm zoltan\@alpha.lehman.cuny.edu}
\endaddress
\abstract This article  
concludes the comprehensive
study started in \cite{Sz5}, where
the first non-trivial isospectral pairs of metrics are constructed 
on balls and  
spheres. These investigations incorporate 4 different cases since these
balls and spheres are considered both on
2-step nilpotent Lie groups and on
their solvable extensions. 
In \cite{Sz5} the considerations are 
completely concluded in the ball-case and in the nilpotent-case. 
The other cases were
mostly outlined.
In this paper
the isospectrality theorems are completely established on spheres.  
Also the important details required about the solvable extensions 
are concluded in this paper. 

A new so called 
{\it anticommutator technique} is developed 
for these constructions. This 
tool is completely different from the other methods
applied on the field so far. It brought a wide range of new
isospectrality examples. 
Those constructed on
geodesic spheres of certain harmonic manifolds 
are particularly striking.
One of
these spheres is homogeneous (since it 
is the geodesic sphere of a 2-point
homogeneous space) while the other spheres, although isospectral to the
previous one, are not even locally homogeneous. This shows that how
little information is encoded about the geometry (for instance, about
the isometries) in the spectrum of Laplacian acting on
functions.
\endabstract

\endtopmatter

\document
Research in spectral geometry started out in the early 60's. This field
might as well be
called audible versus non-audible geometry. This
designation much more readily suggests 
the fundamental question of the field: 
To what 
extent is
the geometry of compact Riemann manifolds encoded in
the spectrum of the Laplacian acting on functions?

It started booming in the 80's, however, 
all the isospectral metrics
constructed 
until the early 90's
had the same local geometry and they differed from each 
other only by their global invariants, such as fundamental 
groups.

Then, in 1993, the first examples of isospectral pairs of metrics 
with different
local geometries were constructed both on closed manifolds 
\cite{G1} and on manifolds with boundaries 
\cite{Sz3,4}.
Gordon established her examples on closed nil-manifolds 
(which were diffeomorphic to tori) 
while this author performed his constructions on topologically 
trivial principal torus bundles over
balls, i. e., on $B^m\times T^3$. 
The boundaries of the latter manifolds are the torus bundles  
$S^{m-1}\times T^3$. The 
isospectrality proofs
are completely different in these two cases. 
On manifolds with boundaries the proof
was based on an explicit computation of the spectrum.The main tool
in these computations was the Fourier-Weierstrass decomposition of the 
$L^2$-function space on the torus fibres $T_p^3$. 

The results of this author were first announced during the San Antonio
AMS Meeting, which was held 
January 13-16, 1993 (cf. Notices of AMS, Dec. 1992,
vol 39(10), p. 1245) and, thereafter, in several seminar talks
given at the University
of Pennsylvania, Rutgers University 
and at the Spectral
Geometry Festival held at MSRI(Berkeley),
in November, 1993. It was circulated in preprint form, 
however, it was published much later
\cite{Sz4}. The later publication includes new materials, such as
establishing the isospectrality theorem
on the boundaries $S^{m-1}\times T^3$ of the
considered manifolds as well.

The author's construction strongly related to the 
Lichnerowicz Conjecture (1946)
concerning harmonic manifolds.
This connection is strongly present also in this paper since
the striking examples offered below 
also relate to the Conjecture.

A {\it Riemann manifold} is said to be {\it harmonic} if its harmonic
functions yield the classical mean value theorem. One can easily
establish this harmonicity on two-point homogeneous manifolds.
The Conjecture claims this statement also in the opposite
direction: The harmonic manifolds are exactly the two-point homogeneous
spaces.

The Conjecture was established on compact, simply connected manifolds
by this author \cite{Sz1}, in 1990. Then, in 1991, Damek and Ricci
\cite{DR} found infinitely many counter examples for the Conjecture 
in the non-compact case by proving that the natural left-invariant
metrics on the solvable extensions of Heisenberg type groups
are harmonic. The Heisenberg type groups are particular 2-step
nilpotent groups attached to Clifford
modules (i .e., to representations of Clifford algebras) \cite{Ka}.
Among them are the groups $H^{(a,b)}_3$ defined by imaginary 
quaternionic numbers (cf. (2.13) and below).

In constructing the 
isospectrality examples described in \cite{Sz3,4}, 
the center $\bold R^3$ of these
groups was factorized by a full lattice $\Gamma$ 
to obtain the torus
$T^3=\Gamma\backslash\bold R^3$ and the torus bundle 
$\bold R^{4(a+b)}\times T^3=\Gamma\backslash H_3^{(a,b)}$. 
Then this torus bundle was restricted onto a ball
$B\subset \bold R^{4(a+b)}$ and both
the Dirichlet and Neumann spectrum of the bundle
$B\times T^3$ 
(topological product)
was computed. It turned out
that both spectra depended only on the value $(a+b)$, 
proving the desired
isospectrality theorem for the ball$\times$torus-type
domains of the metric groups $H^{(a,b)}_3$
having the same $(a+b)$.
 
Gordon and Wilson \cite{GW3} generalized the isospectrality result
of \cite{Sz3,4} to the the ball$\times$torus-type domains 
of general 2-step nilpotent Lie groups. Such a Lie group is uniquely
determined by picking a linear space, $E$, of skew endomorphisms acting
on a Euclidean space $\bold R^m$ (cf. formula (0.1)). Two endomorphism
spaces are said to be {\it spectrally equivalent} 
if there exists an orthogonal transformation between them 
which corresponds
isospectral (conjugate) endomorphisms to each other.
(This basic concept of
this field was introduced
in \cite{GW3}. Note that the endomorphism
spaces belonging to the Heisenberg type groups $H^{(a,b)}_3$ satisfying
$(a+b)=constant$ are spectrally eqivalent.) 

By the first main theorem of \cite{GW3},
{\it the corresponding ball$\times$torus domains
are both Dirichlet and Neumann isospectral on 2-step nilpotent
Lie groups which are defined by spectrally 
equivalent endomorphism spaces.} Then
this general theorem is used for  
constructing
continuous families of isospectral metrics on $B^m\times T^2$
such that the distinct family members have different local
geometries. 

It turned out too that these metrics induce non-trivial
isospectral metrics also on the boundaries, 
$S^{m-1}\times T$, 
of these manifolds. 
This statement was independently established
both with respect to the \cite{GW3}-examples (in
\cite{GGSWW}) and the \cite{Sz3}-examples (in
\cite{Sz4}). Each of these examples have
their own interesting new features. Article \cite{GGSWW} provides the 
first continuous families of isospectral metrics on closed manifolds 
such that the distinct
family members have different local geometries. In \cite{Sz3} 
one has only a discrete family $g^{(a,b)}_3$ 
of isospectral metrics on $S^{4(a+b)}\times T^3$
(such a family is defined by the constant $a+b$).
The surprising new feature is that the metric $g^{(a+b,0)}_3$ is
homogeneous while the metrics $g^{(a,b)}_3$ satisfying $ab\not =0$
are locally inhomogeneous.  

At this point of the development no non-trivial isospectral metrics
constructed on simply connected manifolds were known in the literature.
The first such examples were constructed by Schueth
\cite{Sch1}. The main idea of her construction is the following:
She enlarges the torus $T^2$ of the torus bundle
$S^{m-1}\times T^2$
considered in \cite{GGSWW} 
into a compact simply connected Lie group $S$ such that $T^2$ is a
maximal torus in $S$.
Then the isospectral metrics are constructed on
the enlarged manifold 
$S^{m-1}\times S$.
Also this enlarged manifold is a $T^2$-bundle with respect to the left
action of $T^2$ on the second factor. The original bundle
$S^{m-1}\times T^2$ is
a sub-bundle in this enlarged bundle. 
Then the one parametric families of isospectral metrics introduced in 
\cite{GGSWW} 
on manifolds $S^{m-1}\times T^2$ 
are 
extended 
such that they provide
isospectral metrics also 
on the enlarged manifold.
In special cases she gets examples on
the product of spheres. The metrics with the lowest dimension are
constructed on $S^4\times S^3\times S^3$. 

In \cite{Sch2} 
this technique is reformulated in a more
general form such that certain
principal torus bundles are
considered with a fixed metric on the base space and 
with the natural flat metric 
on the torus $T$. (Important basic concepts of this general theory
are abstracted from works \cite{G2;GW3}.)
The
isospectral metrics are constructed on the total space
such that they have the following three properties: (1) 
The elements of the structure group $T$ act as isometries. (2)
The torus-fibers
have the prescribed natural flat metric. (3) The
projection onto the base space is a Riemannian
submersion. 

One can define such a Riemannian metric just
by choosing a connection on this principal torus bundle for 
defining the orthogonal
complement to the torus fibers. Then the
isospectral metrics with different local geometries 
were found by appropriate 
changing (deforming) of these 
connections. 
This combination of extension- and 
connection-techniques is a key feature of Schueth's constructions,
which provided new surprising examples including
isospectral pairs of metrics with different local geometries with the
lowest known dimension on $S^2\times T^2$.

Let us mention that in each of the papers 
\cite{G1;2;GW3;GGSWW;Sch1} the general torus bundles involved to the
constructions have total geodesic torus fibres. This assumption is
not used in establishing the isospectrality theorem on the special
torus bundle considered \cite{Sz3,4}. This assumption is removed
and the torus bundle technique is formulated in a very general form
in \cite{GSz}. Though this form of the general isospectrality
theorem opens up new directions,
yet examples
constructed on balls or on spheres were still out of touch by this
technique,
since no ball or sphere can be
considered as the total space of a torus bundle, where $dim(T)\geq 2$. 

The first
examples of isospectral metrics on balls and spheres 
have been constructed
most recently by this author \cite{Sz4} and, very soon thereafter, 
by Gordon \cite{G3}
independently.
The techniques applied in these two constructions are completely
different, providing 
completely different examples of isospectral metrics.
Actually none of these examples can be constructed by the technique used
for constructing the other type of examples. 

First we describe Gordon's examples. The crucial new idea in
her construction is a generalization
of the torus bundle technique such that, 
instead of a principal torus bundle,
just a torus action is considered 
which is not required to be free anymore.
Yet this generalization is benefited by 
the results and methods of the bundle 
technique (for instance, by the Fourier-Weierstrass 
decomposition of function spaces
on the torus fibres for establishing the isospectrality theorem) 
since they are still applicable
on the everywhere dense open subset covered by the maximal dimensional
principal torus-orbits. This idea 
really gives the chance for constructing
appropriate isospectral metrics on balls and sphere, 
since these
manifolds admit such non-free torus actions. 

In her construction Gordon uses the 
metrics defined on $B\times T^l$ resp. 
$S\times T^l$ introduced 
in \cite{GW3} resp. \cite{GGSWW}. 
First, she represents the torus $T^l=\bold Z^l\backslash\bold R^l$
in $SO(2l)$ by using the natural identification $T^l=\times_l SO(2)$. By
this representation she gets an enlarged bundle with the
base space $\bold v=\bold R^k$ 
and with the total space $\bold R^{k+2l}$ such that the torus is 
non-freely acting
on the total space. Then a
metric is 
defined on the total space. This metric
inherits the Euclidean metric of
the torus orbits and its projection onto the base space is the original
Euclidean metric. Therefore, only the horizontal subspaces 
(which are perpendicular to the orbits)
should be defined. They are introduced
by the alternating bilinear form
$B:\bold R^k\times\bold R^k\to\bold R^l$, 
where $\< B(X,Y),Z\> =\< J_Z(X),Y\>$. 
Her final conclusion is as follows:

{\it If the one parametric family $g_t$, considered in the
first step on the manifolds $B\times T^l$, or, 
on $S\times T^l$, consists
of isospectral metrics then also 
the above constructed metrics $\wt{g}_t$ are isospectral
on the Euclidean balls and spheres of the total
space $\bold R^{k+2l}$.}

This construction provides locally inhomogeneous metrics 
because the torus actions involved have degenerated orbits.
In the concrete examples, since the metrics $g_t$
constructed in \cite{GW3} and \cite{GGSWW} are used, the 
torus $T$ is 2-dimensional.
In an other theorem Gordon 
proves that the metrics $\wt{g}_t$ 
can be arbitrary close chosen to the
standard metrics of Euclidean balls and spheres.

\head
\leftline{Constructing by the anticommutator technique}
\endhead

The Lie algebra of a 
2-step nilpotent metric Lie group is 
described by a system 
$\{\bold n=\bold v\oplus\bold z,\< ,\> ,J_Z\}$, where the Euclidean
vector space $\bold n$, with the inner product $\< ,\>$, is decomposed
into the indicated orthogonal direct sum, furthermore, $J_Z$ is a skew
endomorphism acting on $\bold v$ for all $Z\in\bold z$ such that the
map $J:\bold z\to End(\bold v)$ is linear and one to one.
The linear space of endomorphisms $J_Z$ is denoted by $J_{\bold z}$.
Then the nilpotent Lie algebra with the center $\bold z$ is defined by
$$
\<[X,Y],Z\>=\<J_Z(X),Y\>\, ;\, 
[X,Y]=\sum_\alpha \<J_{Z_{\alpha}}(X),Y\>Z_\alpha 
\tag 0.1
$$
where $X,Y\in\bold v\, ;\, Z\in\bold z$ and $\{Z_1,\dots ,Z_l\}$ is an
orthonormal basis on $\bold z$.
 
Note that such a Lie algebra is uniquely determined 
by a linear 
space, $J_{\bold z}$, of skew
endomorphisms acting on a Euclidean vector space $\bold v$.
The natural Euclidean norm is defined by 
$||Z||^2=-Tr(J_Z^2)$ on $\bold z$. The constructions below
admit arbitrary other Euclidean 
norm on $\bold z$.

The Lie group defined by this Lie algebra is denoted by $G$. The
Riemann metric, $g$, is defined by the 
left invariant extension of the
above Euclidean inner product introduced on the tangent space 
$T_0(G)=\bold n$ at the origin $0$. The exponential
map identifies the Lie algebra $\bold n$ with the vector space
$\bold v\oplus\bold z$. Explicit formulas 
for geometric objects such as
the invariant vector fields $(\bold{X}_i,\bold{Z}_\alpha )$, Laplacian,
etc. are described in (1.1)-(1.6).

The particular Heisenberg-type nilpotent groups are defined by
special endomorphism spaces satisfying 
$\bold J_Z^2=-|Z|^2id$, for all $Z\in\bold z$
\cite {Ka}. If $l=dim(\bold z)=3mod4$,
then there exist (up to equivalence)
exactly two Heisenberg-type endomorphism spaces, 
$J_l^{(1,0)}$ and
$J_l^{(0,1)}$,
acting irreducibly on $\bold v=\bold{R}^{n_l}$ (see the explanations
at (2.6)). The 
reducible endomorphism spaces can be described by an
appropriate Cartesian product in the form $J_l^{(a,b)}$ (see more 
about this notation below (2.14)). When
quaternionic- resp. Cayley-numbers are used for constructions, 
the corresponding 
endomorphism spaces are denoted by
$J_3^{(a,b)}$ resp.
$J_7^{(a,b)}$.
The family $J_l^{(a,b)}$, defined by fixed values of $l$
and $(a+b)$, consists of spectrally equivalent endomorphism spaces.
 
Any 2-step nilpotent Lie group $N$ extends to a solvable group
$SN$ defined on the half space $\bold n\times\bR_+$ (cf. (1.8) and 
(1.9)).
The first spectral investigations 
on these solvable extensions are established 
in \cite{GSz}. 

The ball$\times$torus-type domains,
which were sketchily introduced above, 
are defined by the factor manifold $\Ga_Z\setminus\bold n$,
where $\Ga_Z$ is a full lattice on the Z-space 
$\bold z$ such that this principal torus bundle is considered
over a Euclidean ball
$B_\de$ of radius $\de$ around 
the origin of the X-space $\bold v$. 
The boundary of this manifold is the principal torus 
bundle $(S_\de,T)$. 

The main tool in proving the isospectrality theorem
on such domains is the Fourier-Weierstrass
decomposition $W=\oplus_\alpha W^\alpha$ of the $L^2$ function
space on the group $G$, where, in the nilpotent case, 
the $W^\alpha$ is spanned by the
functions of the form $F(X,Z)=f(X)e^{-2\pi\sqrt{-1}\<Z_\alpha ,Z\>}$.
It turns out that each $W^\alpha$ is invariant under the action of the
Laplacian, $(\Delta_GF)(X,Z)=\square_{\alpha}(f)(X)
e^{-2\pi\sqrt{-1}\<Z_\alpha ,Z\>}$, such that 
$\square_{\alpha}$ depends,
besides some universal terms and $\Delta_X$, only on 
$J_{Z_\alpha}$ and it does not depend on the other endomorphisms.
Since $J_{Z_\alpha}$ and
$J_{Z_\alpha^\prime}$ are isospectral, one can intertwine the Laplacian
on the subspaces $W^\alpha$ separately by the orthogonal transformation
conjugating  
$J_{Z_\alpha}$ to
$J_{Z_\alpha^\prime}$.
This tool extends not only to the general ball$\times$torus-cases
considered in \cite{GSz} but also to the torus-action-cases considered
in \cite{G3}. 

The simplicity of the isospectrality 
proofs by the above described Z-Fourier transform is due
to the fact that, on an invariant subspace $W^\alpha$,
one should deal only with one endomorphism, $J_\alpha$, while the
others are eliminated. 

New, so called {\it ball-type domains} were introduced in \cite{Sz5}
whose spectral investigation has no prior history. These domains are 
diffeomorphic to Euclidean balls whose smooth boundaries are described
as level sets by equations 
of the form $\varphi (|X|,Z)=0$ resp. 
$\varphi (|X|,Z,t)=0$, according to 
the nilpotent resp.
solvable cases. The boundaries of these
domains are diffeomorphic to Euclidean spheres which are called 
sphere-type manifolds, or, sphere-type
hypersurfaces. 

The technique of the Z-Fourier transform breaks down
on these domains and hypersurfaces, since the functions 
gotten by this transform do not satisfy the required 
boundary conditions. The Fourier-Weierstrass 
decomposition does not apply on the sphere-type hypersurfaces 
either. 
The difficulties in proving the isospectrality 
on these domains originate from 
the fact that no such Laplacian-invariant 
decomposition of the
corresponding $L^2$ function spaces
is known which keeps, 
on an invariant subspace, 
only one of
the endomorphisms active 
while it gets rid of the
other endomorphisms. 
The isospectrality proofs on these manifolds require a new
technique whose brief description is as follows.

Let us mention first that a wide range of spectrally equivalent 
endomorphism spaces were introduced in \cite{Sz5} 
by means of the so called $\sigma$-{\it deformations}. 
These deformations are defined
by an involutive orthogonal transformation $\sigma$ 
on $\bold v$ which commutes
with all of the endomorphisms from $J_{\bold z}$. The $\sigma$-deformed
endomorphism space, $J^\sigma_{\bold z}$, 
consists of endomorphisms of the form $\sigma J_Z$. This new 
endomorphism space is clearly 
spectrally equivalent to the old one. 
Note that no restriction on
$dim(\bold z)$ is imposed in this case. 
These deformations are of discrete
type, however, which can be considered as the generalizations of
deformations considered on the endomorphism spaces $J^{(a,b)}_l$ in
\cite{Sz3,4}. These deformations provide isospectral
metrics on the ball$\times$torus-type domains by the 
Gordon-Wilson theorem.

The new so-called {\it anticommutator technique}, developed for
establishing the spectral investigations on ball- and sphere-type 
manifolds, does not apply for all the $\sigma$-deformations.
We can accomplish the isospectrality theorems by this
technique only for those
particular endomorphism spaces which include non-trivial
{\it anticommutators}. 

A non-degenerated
endomorphism $A\in J_{\bold z}$ is an anticommutator
if and only if $A\circ B=-B\circ A$ holds for all 
$B\in J_{\bold A^\perp}$. 
If an endomorphism space $J_{\bold z}$ contains
an anticommutator $A$, then, by the 
Reduction Theorem 4.1 of \cite{Sz5}, 
a $\sigma$-deformation is equivalent to
the simpler deformation where
one performs $\si$-deformation only on
the anticommutator $A$.
I. e., only $A$ is switched to $A^\sigma =\sigma\circ A$ and one
keeps the orthogonal 
complement $J_{\bold A^\perp}$ 
unchanged.
In \cite{Sz5} and in this paper 
the isospectrality theorems are established for such, so called,
$\sigma_A$-deformations. 

The constructions concern four different
cases, since we perform them on the ball- and sphere-type
domains both of 2-step nilpotent Lie groups and their solvable 
extensions. 
The details are shared
between these two papers. 
Roughly speaking, the proofs are completely
established in \cite{Sz5} on the ball-type domains and all the 
technical details are complete 
on 2-step nilpotent
groups. Though the other cases were outlined in some extend,
the important details concerning the sphere-type domains and
the solvable extensions are left to this paper.  

This paper starts by a review about the solvable extensions of
2-step nilpotent groups. Then, 
in Proposition 2.1, we 
describe all 
the endomorphism spaces having an anticommutator $A$
(alias $ESW_A$'s) in a representation theorem, where the 
Pauli matrices play very crucial role. The basic examples
of $ESW_A$'s are
the endomorphism spaces $J^{(a,b)}_l$ belonging to Clifford modules.
In this case each endomorphism is an anticommutator. 
The representation theorem describes a great abundance of 
other examples.

In Section 2 
the so called $unit_{endo}$-deformations are introduced just  
by choosing two different unit anticommutators $A_0$ and $B_0$ to
a fixed endomorphism space $\bold F$ (the corresponding $ESW_A$'s are
$\bold RA_0\oplus\bold F$  and
$\bold RB_0\oplus\bold F$). Also these deformations can be used
for isospectrality constructions. By clarifying a  
strong connection between $unit_{endo}$- and 
$\sigma_A$-deformations (cf. Theorem 2.2) we point out
that the anticommutator technique
is a discrete isospectral construction technique. In fact, we 
prove that 
continuous 
$unit_{endo}$-deformations provide conjugate
$ESW_A$'s and therefore the corresponding metrics are isometric.

The main isospectrality
theorem is
stated in the following form in this paper.  
 
\proclaim{Main Theorems 3.2 and 3.4} Let 
$J_{\bold z}=J_A\oplus J_{A^\perp}$ and
$J_{\bold z^\prime}=J_{A^\prime}\oplus J_{A^{\prime\perp}}$
be endomorphism
spaces 
acting on the same space $\bold v$ such that  
$J_{A^\perp}=J_{A^{\prime\perp}}$, furthermore, the 
anticommutators $J_A$ and $J_{A^\prime}$ are either unit endomorphisms
(i. e., $A^2=(A^\prime )^2=-id)$ or they are $\sigma$-equivalent.
Then the map 
$\pa\kappa =T^\prime\circ\pa\kappa^* T^{-1}$ intertwines the 
corresponding Laplacians on the sphere-type boundary $\pa B$
of any ball-type domain, both on the nilpotent groups 
$N_J$ and $N_{J^\prime}$ 
and/or on their solvable extensions $SN_J$ and $SN_{J^\prime}$.
Therefore the corresponding metrics are isospectral on
these sphere-type manifolds.
\endproclaim

In \cite{Sz5}, the corresponding theorem is established  only 
for balls and for $\sigma_A$ deformations. The investigations on
spheres are just outlined and even these sketchy details 
concentrate mostly on the striking examples.

The constructions of the intertwining operators 
$\kappa$ and $\pa\kappa$
require an appropriate decomposition of the function
spaces. This decomposition is, however,  
completely different from the Fourier-Weierstrass 
decomposition applied in the torus-bundle cases since
this decomposition is performed 
on the $L^2$-function space of the X-space. The details are as follows. 

The crucial terms in the Laplacian acting on
the X-space are the Euclidean Laplacian $\Delta_X$ and the operators
$D_A\bullet ,D_F\bullet$ derived from the endomorphisms (cf.
(1.5), (1.12), (3,7), (3.33)). The latter
operators commute with $\Delta_X$.  
In the first step only the operators
$\Delta_X$ and $D_A\bullet$ are considered and
a common eigensubspace decomposition of the corresponding $L^2$
function space is established.
This decomposition results a refined
decomposition of the spherical harmonics on the spheres of the X-space.
Then the operators $\kappa ,\pa\kappa$ are defined such that they
preserve this decomposition. Though one can not get rid of the other
operators $D_F\bullet$ by this decomposition,
the anticommutativity of $A$ by the perpendicular 
endomorphisms $F$ ensures that also the terms containing the
operators $D_F\bullet$ in the Laplacian are intertwined by $\kappa$
and $\pa\kappa$.

By proving also
the appropriate non-isometry theorems, 
these examples provide a wide range of
isospectral pairs of metrics constructed on spheres
with different local geometries. These non-isometry proofs are 
achieved by an independent {\it Extension Theorem} asserting that
an isometry between two sphere-type manifolds extends to
an isometry between the ambient
manifolds. (In order to avoid an 
even more complicated proof, the theorem is established for 
sphere-type manifolds described by equations of the form
$\varphi (|X|,|Z|)=0$ resp. 
$\varphi (|X|,|Z|,t)=0$. 
It is highly probable that one can establish
this extension in the most general cases by 
extending the method applied here.)
This theorem traces back the problem of non-isometry to the 
ambient manifolds, where the non-isometry was thoroughly 
investigated in \cite{Sz5}. The extension can be used also for
determining the isometries of a sphere-type manifold by the isometries
acting on the ambient manifold.

The abundance of the isospectral pairs of metrics constructed by
the anticommutator technique on spheres with different local
geometries is exhibited in {\it Cornucopia Theorem 4.9}, which is the
combination of the isospectrality theorems and of the non-isometry
theorems.

These isospectral pairs include
the so called {\it striking examples} construct\-ed on the
geodesic spheres of the solvable groups $SH^{(a,b)}_3$.
(These examples are outlined in \cite{Sz5} with fairly complete details,
yet some of these details are left to this paper.)
These spheres are homogeneous on the 2-point homogeneous space
$SH^{(a+b,0)}_3$  
while the other spheres on  
$SH^{(a,b)}_3$  
are locally inhomogeneous.
These examples
demonstrate the surprising fact that no information 
about the isometries  
is encoded
in the spectrum of Laplacian acting on functions.

\head
\centerline{\bf 1. Two-step nilpotent Lie algebras and 
their solvable extensions}
\endhead

A metric 2-step nilpotent Lie algebra is described by the system  
$$
\frak n=
\big\{
\bold n=\bv\oplus\bz,\<,\>,J_Z\big\},
\tag 1.1
$$
where $\<,\>$ is an inner product defined on the algebra 
$\bold n$, 
the space $\bz =[\bold n,\bold n]$ is the center of $\bold n$,
furthermore 
$\bv$ is the orthogonal complement to $\bz$.
The map
$J :\bz\to SkewEndo(\bv )$ is defined by $\<J_Z(X),Y\>=\<Z,[X,Y]\>$.
 
The vector spaces $\bv$ and $\bz$ are called 
X-space and Z-space respectively.

Such 
a Lie algebra is well defined by the 
endomorphisms $J_Z$. 
The linear space
of these endomorphisms is denoted by $J_\bz$. For a fixed X-vector 
$X\in\bold v$,
the subspace 
spanned by the X-vectors
$J_Z(X)$ (for all $Z\in\bz$) is denoted by
$J_\bz(X)$.

Consider the orthonormal bases $\big\{E_1;\dots;E_k\}$ and 
$\big\{e_1;\dots;e_l\big\}$ on the X- and the Z-space respectively.
The corresponding coordinate systems defined by these bases are 
denoted by
$\big\{x^1;\dots;x^k\big\}$ and $\big\{z^1;\dots;z^l\big\}$.
According to \cite {Sz5} the left-invariant
extensions of the vectors $E_i;e_{\alpha}$ are the 
vector fields
$$
\gathered
\bold
X_i=\pa_i + \frac 1 {2} \sum_{\alpha =1}^l
\<[X,E_i],e_{\alpha} \>  \pa_\alpha = 
\\
=\pa_i + \frac 1 {2} \sum_{{\alpha} =1}^l \<
J_\alpha\big(X\big),E_i\>
\pa_\alpha\quad ; \quad
\bold Z_{\alpha}=\pa_\alpha,
\endgathered
\tag 1.2
$$
where $\pa_i=\pa /\pa x^i$, $\pa_\alpha=\pa/\pa z^\alpha$
and $J_\alpha = J_{e_\alpha}$.

The covariant derivative acting on 
invariant vector fields is described as follows.
$$
\na_XX^*=\frac 1 {2} [X,X^*]\, ;\,\na_XZ=\na_ZX=-\frac 1 {2}
J_Z\big (X\big)\, ;\,\na_ZZ^*=0.
\tag 1.3
$$

The Laplacian, $\Delta$, acting on functions can 
be explicitly established by substituting (1.2) and (1.3) into
the following well known formula
$$
\De=\sum_{i=1}^k\big(\bold X_i^2-\na_{\bold X_i}\bold X_i\big)
+\sum_{\al=1}^l\big (\bold Z_{\al}^2-\na_{\bold Z_{\al}} 
\bold Z_{\al}\big ).
\tag 1.4
$$
 Then we obtain 
$$
\De=\De_X+\De_Z+\frac 1 {4} \sum_{\al,\beta =1}^l \<J_\al
\big (X\big),J_\beta\big (X\big)\>
 \pa_{\alpha\beta}^2
+\sum_{\al=1}^l\pa_\alpha D_\al \bullet,
\tag 1.5
$$
where $D_\al \bullet$ means differentiation (directional derivative)
with respect to the vector field
$$
D_\al : X \to J_{\al}\big (X\big )
\tag 1.6
$$
tangent to the X-space, furthermore 
$\pa_{\alpha\beta}=\pa^2/\pa z^\alpha
\pa z^\beta$.\par Some other basic objects (such as 
Riemannian curvature, Ricci curvature,  
$d$- and $\delta$-operator
acting on forms) are also explicitly
established in \cite {Sz5}. Finally, we mention a theorem describing 
the isometries on 2-step nilpotent Lie groups.

\proclaim {Proposition 1.1 (\cite {K} \cite {E} \cite{GW3} \cite {W}) }
The 2-step nilpotent metric Lie groups $\big(N,g\big)$ and $\big(N',g'
\big)$ are isometric if and only if 
there exist orthogonal transformations 
$A : \bold v \to
\bold {v^\prime }$ and $C:\bold z \to \bold {z^\prime }$ such that
$$
A J_Z A^{-1} = J^\prime_{C(Z)}
\tag 1.7
$$
holds for any $Z \in \bold z$.
\endproclaim

Any 2-step nilpotent Lie group, $N$, extends to a
solvable group, 
$SN$, defined on the half-space $\bold n\times \bR_+$ with
multiplication given by
$$
(X,Z,t)(X^\prime ,Z^\prime ,t^\prime )
=(X+t^{\frac 1 {2}}X^\prime ,Z+tZ^\prime +\frac 1 {2}t^{\frac 1 {2}}
[X,X^\prime ],tt^\prime ).
\tag 1.8
$$
This formula provides the 
multiplication also
on the nilpotent group $N$, since the latter is a subgroup 
determined by $t=1$.
 
The Lie algebra of this
solvable group is $\bold s=\bold n\oplus \bold t$. 
The Lie 
bracket is completely determined by the formulas 
$$
[\pa_t,X]=\frac 1 {2}X\quad ;\quad [\pa_t,Z]=Z\quad;\quad 
[\bold n,\bold n]_{/SN}
=[\bold n,\bold n]_{/N},
\tag 1.9
$$
where $X \in \bold v$ and $Z \in \bold z$.

In \cite{GSz},
a scaled
inner product $\<\, ,\,\>_c$ with scaling factor $c > 0$ is introduced
on $\bold s$ defined 
by the rescaling 
$|\pa_t|=c^{-1}$ and by
keeping the inner product on
$\bold n$ as well as keeping the relation $\pa_t\perp \bold n$. 
The left invariant
extension of this inner product is denoted by $g_c$.

The left-invariant extensions $\bold Y_i, \bold V_\al ,\bold T$ of the
unit vectors 
$$
E_i=\pa_i\quad ,\quad e_\al =\pa_\al \quad ,\quad\epsilon =c\pa_t
$$
at the
origin are
$$
\bold Y_i=t^{\frac 1{2}}\bold X_i\quad ;\quad \bold V_\al =t\bold Z_\al 
\quad ;\quad \bold T=ct\pa_t ,
\tag 1.10
$$
where $\bold X_i$ and $\bold Z_\al$ are the 
invariant vector fields
on N (cf. (1.2)). 

One can establish these latter
formulas by the following standard computations. Consider the 
vectors $\pa_i ,\pa_{\alpha}$ and $\pa_t$ at the origin $(0,0,1)$ such
that they are the
tangent vectors of the curves $c_A (s)=(0,0,1)+s\pa_A$,
where $A=i,\al ,t$. Then
transform these curves to an arbitrary point by left 
multiplications described in (1.8). Then the tangent of the transformed
curve gives the desired left invariant vector at an arbitrary point.  

The covariant derivative can be computed by the well known formula
$$
\<\na_PQ,R\> =\frac 1{2}\{\<P,[R,Q]\>+\<Q,[R,P]\>+\<[P,Q],R\>\},
$$
where $P,Q,R$ are invariant vector fields.
Then we get
$$
\gathered
\na_{X+Z}(X^*+Z^*)=\na^N_{X+Z}(X^*+Z^*)+c(\frac 1{2}
\<X,X^*\>+\<Z,Z^*\>)\bold T ;\\
\na_X\bold T=-\frac c{2}X\quad ;
\quad\na_Z\bold T=-cZ \quad ;\quad
\na_TX=\na_TZ=\na_TT=0,
\endgathered
\tag 1.11
$$
where $\na^N$ is the covariant derivative on $N$
(cf. (1.3)) and $$X,X^*\in\bold v;Z,Z^*\in\bold z;T\in\bold t.$$

The Laplacian on these solvable groups can be established by the same
computation performed on $N$. Then we get  
$$
\gathered
\Delta=t\Delta_X+t^{\frac 1 {2}}\Delta_Z+
\frac 1 {4}t\sum_{\al ;\beta =1}^l
\<J_\al (X),J_\beta (X)\>\pa^2_{\al \beta} \\
+t\sum_{\al =1}^l\pa_\al D_\al\bullet+c^2t^2\pa^2_t+
c^2(1-\frac k {2}-l)t\pa_t.
\endgathered
\tag 1.12
$$

Also the Riemannian curvature can be computed 
straightforwardly
such that formulas (1.11) are substituted into the standard
formula of the Riemannian curvature. Then we get
$$
\gathered
R_c(X^*\wedge X)=R(X^*\wedge X)- \frac c{2}[X^*,X]\wedge\bold T+
\frac {c^2}{4}X^*\wedge X ;
\\
R_c(X\wedge Z)=R(X\wedge Z)-\frac c{4}J_Z(X)\wedge\bold T+\frac{c^2}{2}
X\wedge Z ;
\\ 
R_c(Z^*\wedge Z)=R(Z^*\wedge Z)+c^2Z^*\wedge Z;
\\
R_c((X+Z),\bold T)(.)=c\na_{\frac 1{2}X+Z} (.); \\
R_c((X+Z)\wedge\bold T)={1\over2}c(\sum_\al J_\al(X)\wedge e_\al 
-J^*_Z)+c^2({1\over4}X+Z)\wedge\bold T ,
\endgathered
\tag 1.13
$$
where $J^*_Z$ is the 2-vector dual to the 
2-form $\<J_Z(X_1),X_2\>$ and
$R$ is the Riemannian curvature on 
$N$, described by

$$
\gathered
R(X,Y)X^*=\frac 1 {2} J_{[X,Y]}(X^*) - 
\frac 1 {4} J_{[Y,X^*]}
(X) + \frac 1 {4} J_{[X,X^*]}(Y); \\
R(X,Y)Z=-\frac 1 {4} [X,J_Z(Y)]+
\frac 1 {4} [Y,J_Z(X)]\quad ;
\quad R(Z_1,Z_2)Z_3=0; \\
R(X,Z)Y=-\frac 1 {4} [X,J_Z(Y)]\quad ; \quad 
R(X,Z)Z^*=-\frac 1 {4}
J_ZJ_{Z^*}(X); \\
R(Z,Z^*)X=-\frac 1 {4} J_{Z^*}J_Z(X) + \frac 1 {4}J_Z
J_{Z^*}(X),
\endgathered
\tag 1.14
$$
where $X;X^*;Y \in \bold v$ and $Z;Z^*;Z_1;Z_2;Z_3 \in \bold z$
are elements of the Lie algebra
(cf. also \cite{E}).

By introducing $H(X,X^*,Z,Z^*):=\<J_Z(X),J_{Z^*}(X^*)\>$,
for the Ricci curvature
we have
$$
\gathered
Ricc_c(X)=Ricc(X)-c^2(\frac k{4}+\frac l{2})X;\\
Ricc_c(Z)=Ricc(Z)-c^2(\frac k{2}+l)Z\quad ;\quad
Ricc_c(T)=-c^2(\frac k{4}+l)T,
\endgathered
\tag 1.15
$$
where the Ricci tensor $Ricc$ on $N$ is described
by formulas

$$
\gathered
Ricc(X,X^*)=-\frac 1 {2} \sum_{\alpha =1}^lH(X,X^*,e_\alpha ,e_\alpha )=
-\frac 1{2}H_{\bold v}(X,X^*);
\\
Ricc(Z,Z^*)=\frac 1 {4} \sum_{i=1}^k
H(E_i,E_i,Z,Z^*)=\frac 1{4}H_{\bold z}(Z,Z^*)
\endgathered
\tag 1.16
$$
and by
$Ricc(X,Z)=0$.

By (1.9) we get that the subspaces $\bold v,\bold z$ and $\bold t$ are 
eigensubspaces of the Ricci curvature operator and, except for
finite many scaling factors $c$, the eigenvalue on $\bold t$
is different from the other eigenvalues. For these scaling 
factors $c$, an isometry $\al :SN_c\to SN^\prime _c$ maps $\bold T$
to $\bold T^\prime$ and for a fixed $t$, the hyper-surface $(N,t)$
is mapped to the hyper-surface $(N^\prime ,t^\prime)$, where 
$t^\prime$ is
an appropriate
fixed 
parameter. 
By using left-products on $SN^\prime$, we may suppose
that the $\al$ maps the origin $(0,1)$ to the origin of $SN^\prime$ and
therefore $\al (N,1)=(N^\prime ,1)$. By (1.2) and (1.10), the 
restrictions of the metric
tensors $g_c$ and $g^\prime _c$
onto the hyper-surfaces $(N,1)$ and $(N^\prime ,1)$
are nothing but the metric tensors $g$ resp. $g^\prime$ on the
nilpotent groups. Thus the $\al$ defines an isometry between $(N,g)$ and
$(N^\prime ,g^\prime)$ and so we have:

\proclaim {Proposition 1.2} Except for finite many scaling factors $c$,
the solvable extensions $(SN,g_c)$ and $(SN^\prime ,g^\prime _c)$ are
locally isometric if and only if the nilpotent metric groups $(N,g)$ and
$(N^\prime ,g^\prime)$ are locally isometric.
\endproclaim

It should be mentioned that
the above assumption about the scaling factor can be dropped. This 
stronger
theorem is
proved in \cite{GSz, Proposition 2.13}
by a completely different 
(much more elaborate) technique.

We conclude this section by considering the spectrum of the
curvature operator acting as a symmetric endomorphism on the 
2-vectors. 
These considerations can be used for establishing 
the non-isometry proofs. These non-isometry proofs will be
established in
many different ways, however, in order to get a deeper insight 
into the realm of non-audible geometry. Even though the next theorem is an 
interesting contribution to this geometry, the understanding of the main
course of this paper is not disturbed by continuing the study  
by the next section.

Two symmetric operators are said to be {\it isotonal} if the
elements of their spectra are the same but the multiplicities
may be different. This property 
is accomplished  
for the curvature operators of
$\sigma^{(a+b)}$-equivalent nilpotent groups in 
\cite{Sz5,
Proposition 5.4}. Now we establish this statement also on the 
solvable extensions of these groups. The technical definition
of the groups $N^{(a,b)}_J$ and the $\sigma^{(a+b)}$-deformations
can be found both in \cite{Sz5} and at formulas (2.12)-(2.14) of this
paper.

In the nilpotent case we used the
following decomposition, which technique extends also to 
the solvable case. 

First
decompose the X-space of the considered
nilpotent Lie-algebras 
in the form
$\bold v=\bold v^{(a)}\oplus\bold v^{(b)}$
such that the involution $\sigma^{(a,b)}$ acts on
$\bold v^{(a)}=\bR^{na}$ by $id$ and it acts on the subspace
$\bold v^{(b)}=\bR^{nb}$ by $-id$. Then the
subspaces 
$$
\gathered
\bold D=(\bold v^{(a)}\wedge \bold v^{(a)})\oplus
(\bold v^{(b)}\wedge\bold v^{(b)})\oplus
(\bold z\wedge\bold z); \\ \bold F=
\bold v^{(a)}\wedge\bold v^{(b)}\quad ;\quad \bold G=
\bold v\wedge\bold z,
\endgathered
\tag 1.17
$$
in $\bold n\wedge\bold n$,
are invariant under the action
of the curvature operator ${R^{ij}}_{kl}$.
The space $\bold F$ is further decomposed into the mixed boxes
$\bold F_{rs}=\bold v_r^{(a)}\wedge\bold v^{(b)}_s$, where 
$\bold v_r^{(a)}$ is the $r^{th}$ component subspace, $\bold R^n$, in
the Cartesian product $\bold v^{(a)}=\times \bold R^n$. Then one
can prove
that the spectrum on such a mixed box is the same on 
$\sigma^{(a+b)}$-equivalent spaces and, furthermore, 
it is the negative of the 
spectrum on a mixed box $\bold v_p^{(a)}\wedge\bold v_q^{(a)};\bold
v^{(b)}_p\wedge\bold v^{(b)}_q\subseteq\bold D$, where $p\not =q$. 
These latter
2-vectors
span the complement space, $\bold{Dg}^\perp$, to the diagonal
space 
$$
\bold {Dg}=(\sum_r\bold v_r^{(a)}\wedge\bold v_r^{(a)})\oplus
(\sum_r\bold v_r^{(b)}\wedge
\bold v_r^{(b)})\oplus (\bold z\wedge\bold z) .
$$
On the invariant space $\bold {Dg}\oplus\bold G$ one can prove that the
spectra of the considered operators are the same, since they are 
isospectral
to the curvature operator on the group $N^{(a+b,0)}_{\bold z}$.
Therefore, comparing the two spectra, we get that only the
multiplicities of eigenvalues belonging to the mixed boxes
of the invariant
spaces $\bold F$ resp. 
$\bold {Dg}^\perp$ are different, 
while the elements
of the spectra are the same. These multiplicities depend on the number
of the mixed boxes, i. e., on $ab$.
This proves that the curvature operator on 
$N^{(a+b,0)}_{\bold z}$
is subtonal to the operator on
$N^{(a,b)}_{\bold z}$ and the curvature operators on
$N^{(a,b)}_{\bold z}$ and
$N^{(a^\prime ,b^\prime)}_{\bold z}$,
where $a+b=a^\prime +b^\prime$ and $aba^\prime b^\prime\not =0$, are
isotonal.

On the solvable extensions, $SN^{(a,b)}_\bold z$, the corresponding
invariant subspaces are the following ones:
$$
\bold F\quad ,\quad \bold{Dg}^\perp\quad ;\quad \bold G\quad ;\quad
\bold{Dg}\oplus (\bold n\wedge\bold t).
\tag 1.18
$$
First consider the last subspace. 
From (1.13) we get that
the map $\tau$, defined by $\tau=-id$ on the
space $(\bold v^{(b)}\wedge\bold v^{(b)})\oplus(\bold v^{(b)}
\wedge\bold t)$ 
and by $\tau =id$ on the orthogonal complement, intertwines
the curvature operators of the spaces 
$SN_{\bold z}^{(a,b)}$ and
$SN_{\bold z}^{(a+b,0)}$ on this subspace. 
Actually, this statement is true on the
direct sum of $\bold G$ and the subspace listed on the last place
in (1.18). Furthermore, the spectrum 
$\{\nu_i\}$
on a mixed box $\bold F_{pq}$ is the same on $\sigma^{(a+b)}$-equivalent
spaces which can be expressed with the help of the 
corresponding spectrum 
$\{\lambda_i\}$ on the nilpotent group in the form $\nu_i=-Q^2+
\lambda_i$. Then the spectrum on a mixed box of $\bold {Dg}^\perp$
has the form $\{-Q^2-\lambda_i\}$. We get again that only the
multiplicities corresponding to these eigenvalues are different
with respect to the two spectra, since these multiplicities depend
on the number of the mixed boxes (i. e., on $ab$). Thus we have
\proclaim {Proposition 1.3} The curvature operators on the 
$\si$-equivalent metric Lie groups
$SN^{(a,b)}_{\bold z}$ and
$SN^{(a^\prime ,b^\prime )}_{\bold z}$  
with $aba^\prime b^\prime\not =0$ are isotonal.

In many cases they are isotonal 
yet non-isospectral. This is the case,
for instance, on 
the groups $S\bold H^{(a,b)}_3$ 
with the same $a+b$ and $ab\not =0$, where the curvature operators 
are isotonal yet non-isopectral
unless
$(a,b)=(a^\prime ,b^\prime )$ up to an order. 

A general
criteria can be formulated as follows:
The Riemannian curvatures on the spaces 
$SN^{(a,b)}_{\bold z}$ and
$SN^{(a^{\prime},b^{\prime})}_{\bold z}$ with $(a+b)=(a^\prime
+b^{\prime})$ and
$0\not =ab\not =a^\prime b^\prime \not =0$ 
are strictly isotonal 
if and only if the spectrum of the curvature of the corresponding
nilpotent group changes on the
mixed boxes $\bold F_{pq}$ when it is multiplied by $-1$.

The curvature of $SN^{(a+b,0)}_{\bold z}$ is just subtonal
(i. e., 
the tonal spectrum is a proper subset of the other tonal spectrum) to
 the curvatures of the manifolds 
$SN^{(a,b)}_{\bold z}$ with $ab\not =0$.
\endproclaim

\head
\centerline{\bf 2. Endomorphism spaces with anticommutators 
(alias $ESW_A$)} 
\endhead

For accomplishing the isospectrality examples  
a new so called anticommutator technique 
is developed
in \cite{Sz5}. 
A non-degenerated endomorphism $A=J_Z$ is called 
an {\it anticommutator} 
in $J_{\bold z}$  
if
$A\circ B=-B\circ A$ holds for all $B\in J_{Z^\perp}$. I. e., the
endomorphism $A$  
anticommutes with each endomorphisms orthogonal to $A$.

An anticommutator  
satisfying $A^2=-id$ is said to be a
{\it unit anticommutator}. 
Any 
anticommutator can be 
rescaled to a unit anticommutator, since it can be written in the form
$A=S\circ
A_0$, where the symmetric "scaling"
operator $S$ is one of the square-roots of
the operator 
$-A^2$, furthermore, $A_0$ is a unit anticommutator. Then the operator
$S$ is commuting with all elements of the endomorphism space.

The isospectrality examples are accomplished 
by certain deformations performed on $ESW_A$'s. By these deformations
only the $A$ is deformed to 
a new anticommutator 
$A^\prime$ which is 
isospectral (conjugate) to $A$. 
The orthogonal
endomorphisms are kept unchanged (i. e., $A^\perp=A^{\prime\perp}$) and
for a general endomorphism the deformation is defined according to the
direct sum $A\oplus A^\perp$. Such deformations are,
for example, the 
$\sigma_A$-deformations introduced in \cite{Sz5} (seee the definition
also in this paper at (2.12)-(2.14)). An 
other obvious example is 
when both $A$ and $A^\prime$ are unit endomorphisms anticommuting
with the endomorphisms of a given endomorphism space 
$A^\perp =A^{\prime\perp}$. 
We call these deformations
$unit_{endo}$-{\it deformations}. In this paper we consider only these
two sorts of deformations and 
the general isospectral deformations of an 
anticommutator will be studied elsewhere.

A brief outline of this section is as follows. 

First we explicitly
describe all the $ESW_A$'s in a representation theorem, where the
endomorphisms are represented as matrices of Pauli matrices.
(In \cite{Sz5} only particular $ESW_A$'s were constructed to show the
wide range of examples covered by this concept.)

Then the explicit description of the  
$unit_{endo}$-deformations follows.
We also prove
that 
the 
endomorphism spaces $ESW_A$ and $ESW_{A^\prime}$ 
are conjugate  
if $A$ and 
$A^\prime$ can be connected by a continuous curve passing through 
unit anticommutators. 
This statement shows that the anticommutator
technique developed in these papers 
is a discrete construction technique since the corresponding metrics
constructed by continuous $unit_{endo}$ 
deformations
are isometric.  

This section is concluded by describing those $\sigma_A$-
or $unit_{endo}-$deforma\-tions
which provide non-conjugate endomorphism spaces and therefore also
the corresponding metrics are locally non-isometric.

\head
\leftline {The Jordan form of an $ESW_A$}
\endhead

First, we explicitly describe a general
$ESW_A$ by a matrix-representation. Then  
more specific endomorphism spaces such as quaternionic $ESW_A$'s 
(alias $\bold HESW_A$) and Heisenberg-type $ESW_A$'s will be considered.

($\bold A$)
In the following matrix-representation of an
$ESW_A$ the endomorphisms are represented as block-matrices, 
more precisely, they 
are the matrices of 
the following $2\times 2$
matrices (blocks). 
   $$\bold 1=\pmatrix
     1 & 0\\
     0 & 1
     \endpmatrix \,\, ,\,\,
   \bold i=\pmatrix
     0 & 1\\
    -1 & 0
     \endpmatrix \,\, ,\,\,
   \bold j=\pmatrix
    -1 & 0\\
     0 & 1
     \endpmatrix \,\, ,\,\,
   \bold k=\pmatrix
     0 & 1\\
     1 & 0
     \endpmatrix .
\tag 2.1 $$

The matrix product with these matrices are described as follows.
$$
\bold i^2=-\bold 1\, ,\,\bold j^2=\bold k^2=\bold 1\, ,\,\bold i\bold j
=-\bold j\bold i=\bold k\, ,\,\bold k\bold i=-\bold i\bold k=\bold j\,
,\,
\bold k\bold j=-\bold j\bold k=\bold i.
\tag 2.2
$$
The second and the last group of these equations show that (2.1) is
not a 
representation of the quaternionic numbers. Note that the matrices 
$$
\sigma_x=\bold k\quad ,\quad\sigma_y=-\sqrt{-1}\,\bold i\quad ,\quad
\sigma_z=-\bold j
$$ 
are the so called 
Pauli spin matrices. 
 
From the above relations  
the following {\it observation} 
follows immediately:
{\it a $2\times 2$-matrix, $Y$, anticommutes with
$\bold i$ if and only if it has the form 
$Y=y_2\bold j+y_3\bold k$.}

In the following we describe the whole space of skew endomorphisms
anticommuting with a fixed 
skew endomorphism $A$. The endomorphisms
are considered to be represented in matrix form 
such that the matrix of $A$
is a diagonal Jordan matrix. One can establish this 
representation of an $ESW_A$ by an
orthonormal Jordan basis corresponding to the 
anticommutator $A$. 

We consider the eigenvalues of the symmetric endomorphism $A^2$
arranged in the form
$-a^2_1<\dots <-a^2_s\leq 0$.
The corresponding multiplicities are denoted by $m_1,\dots ,m_s$. First
suppose that $A$ is
non-degenerated and therefore 
the main diagonal 
of its Jordan matrix is
built 
up by
$2\times 2$ matrices in the block-form 
$$
(|a_1|\bold i,\dots ,|a_1|\bold i,\dots\dots 
,|a_s|\bold i,\dots ,|a_s|\bold i).
\tag 2.3
$$ 
The $m_c=2n_c$-dimensional eigensubspace belonging to the eigenvalue 
$-a^2_c$
is denoted by 
$\bold B_c^{m_c}$. 

We seek also the anticommuting matrices in the above block-form, i. e.,
we consider them as matrices of $2\times 2$-matrices. Any 
matrix, $F$, 
anticommuting with $A$
leaves the eigensubspaces 
$\bold B_c^{m_c}$ invariant.
Therefore it can be written in the form 
$F=\oplus F_c$, where $F_c$
operates on  
$\bold B_c^{m_c}$.
By the above observation 
we get that $F$ is anticommuting with $A$ if and only if the matrix of
$F_c$, considered as the matrix of $2\times 2$ matrices, has
the block-entries of the form
$F_{cml}=j_{cml}\bold j+k_{cml}\bold k$. 
Since the matrices $\bold j$ and $\bold k$ are 
symmetric, the main diagonal is trivial ($F_{cll}=\bold 0$), 
furthermore,
$j_{cml}=-j_{clm}\, ;\, k_{clm}=-k_{cml}$ hold. I. e., 
an endomorphism $F$
anticommutes with $A$ if and only if the real matrices $j_c$ 
and $k_c$ are skew symmetric.

If $A$ is degenerated, then $a_s=0$ and its action is trivial
on the maximal eigensubspace
$\bold B_s^{m_s}$. In this case the block $F_s$ can be an
arbitrary real skew-matrix.

An {\it irreducible block-decomposition} of an $ESW_A$ is defined
as follows. First we 
decompose
the eigensubspaces 
$\bold B_c^{m_c}$ 
into orthogonal subspaces
$\bold B_{ci}^{m_{ci}}$ such that the endomorphisms leave them invariant
and act on them irreducible. Then we consider a basis whose
elements are in these irreducible spaces. With respect to such a basis, 
all the endomorphisms appear in the form $F_c=\oplus F_{ci}$.
This irreducible decomposition of the X-space is 
the most refined one such that 
the endomorphisms $F$
still can be represented in the form 
$F_{ci}=\{j_{cikl}\bold j +k_{cikl}\bold k\}$. 
The
entry $a_{ci}\bold i$ 
is constant
with  multiplicity $m_{ci}$. 

These statements completely describe the space of skew endomorphisms
anticommuting with $A$. 
If $A$ is non-degenerated, the dimension of this
space is
$\sum_c n_c(n_c-1)$. If $A$ is degenerated, the last 
term in this sum should be changed to $(1/2)m_s(m_s-1)$. A general
$ESW_A$ is an $A$-including subspace of this maximal space.

By summing up we have
\proclaim{Proposition 2.1} Let $\oplus B_c$ be the 
above described Jordan decomposition
of the X-space with respect to an anticommutator $A$ such that $A^2$
has the constant eigenvalue $-a^2_c$ on $B_c$. 
Then all the endomorphisms
from $ESW_A$ leave these Jordan subspaces invariant and,
in case $a_c\not =0$, an
$F\in A^\perp$ can be represented as the matrix of $2\times 2$ matrices
in the form 
$F_c=(F_{cml}=j_{cml}\bold j + k_{cml}\bold k)$, where $j_c$ and $k_c$
are real skew matrices. If $a_c=0$, the matrix representation of 
$F_c$ can be an arbitrary real skew matrix. 

This Jordan decomposition, $ESW_A=\oplus ESW_{A_c}$, can be refined by
decomposing a subspace $B_c$ into irreducible subspaces 
$\bold B_{ci}^{m_{ci}}$. Then any $F_c$ can be represented in the form
$F_c=\oplus F_{ci}$ such that the components, $F_{ci}$, still have the
above described form.  

For a fixed anticommutator $A$, the  
dimension of the maximal $ESW_A$ consisting all the skew endomorphisms
anticommuting with $A$ is
$\sum_{ci} m_{ci}(m_{ci}-1)$. A general $ESW_A$ is an $A$-including
subspace of this maximal endomorphism space. 
\endproclaim 

($\bold B$)
Particular, so called {\it quaternionic
endomorphism spaces with anticommutators} (alias $\bold HESW_A$)
can be introduced
by using matrices with quaternionic entries. 
In this case the X-space is 
the n-dimensional quaternionic vector space identified with
$\bR^{4n}$. 
We suppose that the entries of an
$n\times n$ quaternionic
matrix $A$ act by left side products on the component of
the quaternionic n-vectors.
Such a matrix defines a skew symmetric endomorphism on 
$\bR^{4n}$ if and only if it is a Hermitian skew matrix,
i. e., $a_{ij}=-\overline{a}_{ji}$ holds for the entries.

Notice that in this case 
the entries in the main diagonal are imaginary
quaternions 
furthermore
$A^2$ is a Hermitian symmetric matrix and therefore the
entries in the main diagonal of the matrices $(A^2)^k$ are real
numbers.

A typical example for an $\bold HESW_A$ 
is when $A$ is a diagonal matrix having the 
same imaginary quaternion (say $\bold I$) in the main diagonal
and the anticommuting matrices are symmetric matrices with entries
of the form
$y_2\bold J+y_3\bold K$. If the action of endomorphisms is 
irreducible and
we build up diagonal block matrices by using such blocks, we get
the quaternionic version of the above Proposition 2.1.

Notice that the matrices in a general $ESW_A$
can not be represented as such quaternionic matrices in general. 
In fact, the endomorphisms 
restricted to a subspace $\bold B_{ci}^{m_{ci}}$ 
can be commonly
transformed into
quaternionic matrices if and only if
the multiplicities $m_{ci}$ are the multiples of 4 
($m_{ci}=4k_{ci}$) and
the matrices are matrices of such $4\times 4$ blocks
which are the linear combinations of
matrices of the form
   $$\bold I=\pmatrix
     \bold i &\bold 0\\
     \bold 0 & \bold i
     \endpmatrix \,\, ,\,\,
   \bold J=\pmatrix
    \bold 0 &\bold j\\
    -\bold j & \bold 0
     \endpmatrix \,\, ,\,\,
   \bold K=\pmatrix
     \bold 0 &\bold k\\
    -\bold k &\bold 0
     \endpmatrix .
\tag 2.5 $$
Note that in this quaternionic matrix form, two anticommuting matrices
can be pure diagonal matrices.

($\bold C$)
Other specific {\it endomorphism spaces}
are those {\it where all the
endomorphisms are anticommutators}.

Since on a Heisenberg-type group the equation 
$$ 
J_ZJ_{Z^*}+J_{Z^*}J_Z=-2\<Z,Z^*\>id
$$
holds (cf. (1.4) in \cite {CDKR}, where this statement is proved by
polarizing $J_Z^2=-|Z|^2\,id$),
all the endomorphisms are anticommutators in the endomorphism
space $J_{\bold z}$
of these groups. 

The endomorphism spaces belonging to Heisenberg-type groups are
attached to Clifford modules (which are
representation of Clifford algebras)
\cite{Ka}. Therefore we call them Heisenberg-type, or, Cliffordian
endomorphism spaces. The classification of Clifford modules is well
known, providing classification also for the Cliffordian endomorphism
spaces. Next we briefly summarize some of the main results of 
this theory (cf. \cite{L}).

If $l=dim(J_{\bold z})\not =3(mod 4)$ 
then there exist (up to equivalence) exactly one
irreducible H-type endomorphism space acting on a $\bold R^{n_l}$,
where the dimension $n_l$, depending on $l$, 
is described below. This endomorphism space
is denoted by $J_l^{(1)}$. If $l=3(mod 4)$, then there exist 
(up to equivalence) exactly
two non-equivalent irreducible H-type endomorphism spaces acting on
$\bold R^{n_l}$ which are denoted by 
$J_l^{(1,0)}$ and
$J_l^{(0,1)}$ separately. The values $n_l$ corresponding to
$
l=8p,8p+1,\dots ,8p+7
$
are
$$
n_l=2^{4p}\, ,\, 2^{4p+1}\, , \, 2^{4p+2}\, , \,
2^{4p+2}\, ,
\, 2^{4p+3}\, ,\, 2^{4p+3}\, , \, 2^{4p+3}\, , \,
2^{4p+3}.
\tag 2.6
$$

The reducible Cliffordian endomorphism spaces can be built up by these
irreducible ones. They are denoted by 
$J_l^{(a)}$ resp. $J_l^{(a,b)}$,
corresponding to the definition of
$J_3^{(a,b)}$ and $J_7^{(a,b)}$. (See more
explanation about this notation after formula
(2.14).)

Riehm \cite{R} described these
endomorphism spaces explicitly and used his description to determine
the isometries on Heisenberg-type metric groups.

From our point of view particularly important examples are the 
groups $H^{(a,b)}_3$. The endomorphism space $J^{(a,b)}_3$
of these groups,
defined by appropriate multiplications with imaginary quaternions, 
are thoroughly described in \cite{Sz5}. An other 
interesting case is
$H^{(a,b)}_7$, where 
the imaginary Cayley numbers are
used for constructions. A brief description of this 
Cayley-case
is as follows.

We identify 
the space of imaginary Cayley numbers
with $\bR^7$ 
and we introduce also the maps $\Phi :\bR^7\to \bold H=\bR^4$ and
$\Psi :\bR^7\to\bold H$ defined by
$
(Z_1,\dots,Z_7)\to Z_1\bold i+Z_2\bold j+Z_3\bold k
$
and
$
(Z_1,\dots,Z_7)\to Z_4+ Z_5\bold i+Z_6\bold j+Z_7\bold k
$
respectively. I. e., if we consider 
the natural decomposition $\bold{Ca}=\bold H^2$
on the space $\bold{Ca}$ of
Cayley numbers, 
then the above maps are the corresponding projections onto
the factor spaces. Then the right product $R_Z$ by an imaginary
Cayley number $Z\in\bR^7$ is described by the following formula
$$
R_Z(v_1,v_2)=(v_1\Phi (Z),-v_2\Phi (Z))+(\Psi (Z)v_2,-\overline{\Psi}
(Z)v_1),
$$
where $(v_1,v_2)$ corresponds to the decomposition 
$\bold{Ca}=\bold H^2$. 

The $R_Z$ is a skew symmetric endomorphism satisfying the property
$R^2_Z=-|Z|^2id$ and the whole space of these endomorphisms defines
the Heisenberg type Lie algebra $\frak n^1_7$. 
Then the Lie algebras $\frak n^{(a,b)}_7$ can be similarly defined
than the algebras $\frak n_3^{(a,b)}$. (See more about this notation
below (2.14)).

\head
\leftline{The $unit_{endo}$-deformations
$\delta_{\bold F}:A_0\to B_0$}
\endhead
So far only $\sigma_A$-deformations of an $ESW_A$ were considered.
Seemingly a new type of deformations can be introduced as
follows.

Consider an endomorphism space $\bold F$ spanned by the
orthonormal basis
$\{F^{(1)},\allowmathbreak\dots ,F^{(s)}\}$ and
let $A_0$ and $B_0$ be unit
endomorphisms ($A_0^2=B_0^2=-id$) such that both anticommute with
the elements of $\bold F$. Then the linear map defined by
$A_0\to B_0$ and $F^{(i)}\to F^{(i)}$ between 
$ESW_{A_0}=\bold RA_0\oplus\bold F$ and
$ESW_{B_0}=\bold RB_0\oplus\bold F$ is an orthogonal transformation
corresponding isospectral endomorphisms to each other. The latter
statement immediately follows by
$(F+cA_0)^2=(F+cB_0)^2=F^2-c^2id$.

These transformations are called {\it $unit_{endo}$-deformations} and
are denoted by $\delta_{\bold F}:A_0\to B_0$. Since the isospectrality
theorem extends to these deformations, it is important to compare
them with the $\sigma_A$-deformations. This problem is 
completely answered
by the following theorem.

\proclaim {Theorem 2.2}  
Let $A_0$ resp. $B_0$ unit anticommutators with respect to the 
same system $\bold F=Span\{F^{(1)},\dots ,F^{(s)}\}$. 
Then the orthogonal endomorphism $\sqrt D$,
where the $D$ is derived from $A_0B_0^{-1}$ by 2.7, conjugates $B_0$ to
an anticommutator of $\bold F$ such that it is a $\sigma$-deformation
of $A$. 
Thus any
non-trivial 
$unit_{endo}$-deformation, 
$\delta_{\bold F}:A_0\to B_0$, 
is equivalent to a 
$\sigma_{A_0}$-deformation. 

A continuous family of $unit_{endo}$-deformations
is always trivial. I. e., it is the family of
conjugate endomorphism spaces and therefore the corresponding
metric groups are isometric.
\endproclaim
\demo{Proof}
In this proof we seek after an orthogonal
transformation conjugating $B_0$ to an endomrphism of the form 
$B_0^\prime =\hat S\circ A_0$, where $\hat S$ is a symmetric
endomorphism satisfying $\hat{S}^2=id$, such that the
conjugation fixes, meanwhile, all the endomorphisms from $\bold F$.
Then one can easily establish that the $B_0^\prime$ is the
$\sigma =\hat S$-deformation of $A_0$.
 
The endomorphism
$E=A_0\circ B_0^{-1}=-A_0\circ B_0$ commutes with the $F$'s since they
anticommute both with $A_0$ and $B_0$. 
Decompose $E$ 
into the form
$$
E=A_0\circ B_0^{-1}=-A_0\circ B_0
=S+S^*\circ C=\wt S\circ D,
\tag 2.7
$$
where $S=(1/2)(
A_0\circ B_0^{-1}+
B_0^{-1}\circ A_0)$
is the symmetric part, the endomorphism
$S^*\circ C=(1/2)(
A_0\circ B_0^{-1}-
B_0^{-1}\circ A_0)$
is the skew-symmetric part written in scaled form ($C^2=-id$), and
the orthogonal endomorphism $D$ (commuting with all endomorphisms
 $\{F^{(1)},\dots ,F^{(s)}\}$) is constructed as follows.

Notice that $S$ and $S^*\circ C$ commute and therefore a
common Jordan decomposition
can be established such that the matrix of
$E$ appears as a diagonal matrix of $2\times 2$ matrices
of the form
   $$E_{a}=\pmatrix
     S_a &S^*_a \\
     -S^*_a &S_a
     \endpmatrix =
   \wt S_a\pmatrix
     \cos\al_a &\sin\al_a \\
     -\sin\al_a &\cos\al_a 
     \endpmatrix ,
\tag 2.8$$
where $\wt S_a=(S_a^2+S_a^{*2})^{-1/2}$.
This Jordan decomposition can be described also in the following
more precise form.

The skew endomorphism $[A_0,B_0]$ vanishes exactly
on the subspace $K$, where $A_0$ and $B_0$ commute and therefore 
we should deal only with
these endomorphisms 
on the orthogonal complement $K^\perp$. On this space
the non-degenerated
operators $[A_0,B_0]$ and $B_0$ 
anticommute, generating the quaternionic
numbers 
and both can be represented as diagonal quaternionic matrices
such that 
$C_a=[B_0,A_0]_{0a}=\bold I$
and $B_{0a}=\bold J$ (cf. Lemma 2.4).
Then a $4\times 4$ quaternionic block $E_a$ of $E$ appears in the form
$E_a=\wt S_a(\cos\al_a \, \bold 1+\sin\al_a\,\bold I)$.

On the subspace $K$ (which can be considered
as a complex space with the complex structure $B_0$), the $2\times 2$
Jordan blocks introduced in (2.8) are 
$E_a=S_a\, \bold 1$ and $B_{0a}=\bold i$.

The endomorphism $B_0$ commutes with the symmetric part $D_+$ of $D$
and it is anticommuting with the skew part $D_-$ of $D$. The same 
statement is true with respect to the square root operator $\sqrt D$,
which has the Jordan blocks
 $$\pmatrix
     \cos (\al_a/2) &\sin (\al_a/2) \\
     -\sin (\al_a/2) &\cos (\al_a/2) 
     \endpmatrix .
\tag 2.9$$
Therefore $\sqrt{D}B_0=B_0\sqrt{D}^T=B_0\sqrt{D}^{-1}$ and thus
$$
\sqrt D\circ B_0\circ{\sqrt D}^{-1}=D\circ B_0=\hat S\circ A_0
\tag 2.10
$$
where $\hat S$ is a symmetric, while $\hat SA_0$ is a skew-symmetric
unit endomorphism. Thus $\hat SA_0=A_0\hat S, {\hat S}^2=id$ 
and $B_0^\prime =\hat SA_0$ commutes with $A_0$.

The operator $D$ 
commutes with each of the operators
$\{F^{(1)},\dots ,F^{(s)}\}$, therefore the matrices of the 
$F$'s are symmetric quaternionic block matrices with entries of the
form $f_{ij}\bold I$ such that the blocks $F_{ck}$ corresponding to an
irreducible subspace $B_{ck}$ are included in the blocks
determined by those maximal
eigensubspaces where the values $\wt S_a$ are constant . Therefore also
$\sqrt D$ commutes with these operators. 

Thus the endomorphisms from $\bold F$ anticommute with
$B_0^\prime =\hat SA_0$ and commute with the orthogonal 
transformation $\hat S$. It follows that $B_0^\prime$ is an 
anticommutator with respect to the system $\bold F$ such that it is
the $\sigma =\hat S$-deformation of $A_0$. 

The second part of the theorem obviously follows from the first one.
Thus the proof is concluded.
\qed
\enddemo

The above theorem proves that one can not construct non-trivial
continuous families of isospectral metrics by the 
$unit_{endo}$-deformations.
The following theorem establishes a similar statement corresponding
to the 2-dimensional $ESW_A$'s.

\proclaim{Theorem 2.3}
On a 2-dimensional $ESW_A$ any 
$\sigma_A$-deformation (or unit\-endo-deformation) is
trivial, resulting conjugate endomorphism spaces.
\endproclaim

\demo{Proof}
This theorem is established by the following

\proclaim{Lemma 2.4} 
Let $A$ and $F$ be anticommuting endomorphisms. If both
are non-degenerated, they generate the quaternionic numbers 
and both
can be represented as a diagonal quaternionic matrix  such that 
there are $\bold I$'s on the diagonal  
of $A$ and there are $\bold J$'s on the diagonal of $F$. 
\endproclaim

This lemma easily settles the statement. 

In fact, if the
anticommuting endomorphisms $A$ and $F$ form a basis in the $ESW_A$
such that both are non-degenerated, then represent them in the
above described diagonal quaternionic matrix form. Since the 
irreducible subspaces $B_{ck}$ are nothing 
but the 4-dimensional quaternionic
spaces $\bold H=\bold R^4$, a $\sigma_A$-deformation should operate
such that 
some of the matrices $\bold I$ are switched to $-\bold I$ at some
entries on the diagonal. This operation results the new
endomorphism $A^\prime$.
Let $d_-$ be the set of positions where
these switchings are done. Since $\bold J^{-1}\bold I\bold J=-\bold I$ 
and $\bold J^{-1}\bold J\bold J=\bold J$, 
the quaternionic 
diagonal matrix, having the entry $\bold J$ at a position
listed in $d_-$ and the entry $\bold 1$ at the other positions,
conjugates $A$ to $A^\prime$ while this conjugation fixes $F$.

If one of the endomorphisms, say $F$, is degenerated on a maximal
subspace $K$, then the $A$ leaves this space invariant. If $A$ is
non-degenerated on $K$, then it
defines a complex structure on it. The conjugation by the reflection
in a real subspace takes back $-A_{/K}$ to $A_{/K}$.

The problem of conjugation is trivial on the maximal subspace $L$
where both endomorphisms are degenerated. This proves 
the statement completely.

The proof is concluded by {\it proving Lemma 2.4.}

Represent the non-degenerated endomorphisms $A$ and $F$ 
in the scaled form
$A=S_AA_0\, ;\, F=S_FF_0$. 
Consider also the skew endomorphism $E=AF=S_EE_0$,
where $S_E=S_AS_F$ and $E_0=A_0F_0$. 
It is anticommuting with the endomorphisms $A$ and $F$. Then the 
endomorphisms
$$
A_0=J_{\bold i}\quad ,\quad F_0=J_{\bold j}\quad ,\quad E_0=J_{\bold k}
\tag 2.11
$$
define a quaternionic structure on the X-space and the symmetric
endomorphisms $S_A\, ,\, S_F\, ,\, S_E$ commute with each other as well
as with the 
skew endomorphisms listed above. 

Because of these commutativities, the X-space can be decomposed
into a Cartesian product $\bold v=\oplus\bold H_i$ of pairwise
perpendicular 4-dimensional quaternionic spaces such that all the
above endomorphisms can be represented as diagonal 
quaternionic matrices.
In this matrix form the entries of the matrices
corresponding to the symmetric 
endomorphisms $S$ are real numbers which are
nothing but the eigenvalues of these matrices. From the 
quaternionic representation we get that each of
these eigenvalue-entries has multiplicity 4.
  
This completes the proof both of the Lemma and the Theorem.
\qed
\enddemo

\noindent{\bf Remark 2.5} The isospectrality 
theorem in \cite{Sz5} states
that $\sigma_A$-deform\-ations provide pairs of endomorphisms spaces
such that the ball-type domains with the same radius-function 
are isospectral
on the corresponding nilpotent groups as well as on their solvable 
extensions. 

We would like to modify 
Remark 4.4 of \cite{Sz5},
where the extension of the above isospectrality 
theorem to arbitrary isospectral deformations of an anticommutator
is suggested. 
The spectral investigation of these general deformations
appears to be far more difficult a problem 
than it seen to be earlier. In this
paper we give only a weaker version of this generalization, where
$A$ is supposed to be a unit anticommutator.

This weaker generalization immediately follows from 
Theorem 2.2. 
  
\proclaim{Theorem 2.6} The $ESW_A$-extensions of 
a fixed endomorphism space
$\bold F=Span\allowmathbreak\{F^{(1)},\dots F^{(k)}\}$ by
unit anticommutators $A$ define nilpotent groups (and solvable
extensions) such that for any two of these metric groups
the ball-type
domains with the same X-radius function are isospectral.
\endproclaim

An $ESW_A$-extension of the above fixed set means adding such 
a skew endomorphism $A$ to the system which anticommutes with the
endomorphisms $F^{(i)}$. 

\head
\leftline{$\sigma_A$-deformations providing non-conjugate
$ESW_A$'s}
\endhead

The precise forms of theorems quoted below require 
the precise forms of definitions given for $\sigma-\, ,\,\sigma_A-\, ,\,
\sigma^{(a,b)}$- and 
$\sigma_A^{(a,b)}$-deforma\-tions, performed on an endomorphism space.
These concepts were 
introduced in \cite{Sz5}
as follows.

Let $\sigma$ be an involutive orthogonal transformation commuting with
the endomorphisms of an
$ESW_A=\bold A\oplus\bold A^\perp$, where $\bold A=\bold RA$. 
Then the $\sigma_A$-deformation of
the endomorphism space is defined by deforming $A$ to
$A^\prime =\sigma\circ A$ while keeping the orthogonal endomorphisms
unchanged. The deformation of a general element is defined according to
the direct sum $ESW_A=\bold A\oplus\bold A^\perp$.

These deformations provide spectrally
equivalent endomorphism spaces, since we have
$$
(A^\perp +\sigma A)^2=(A^\perp )^2+A^2=(A^\perp +A)^2.
\tag 2.12
$$
Therefore, there exists an orthogonal transformation between $ESW_A$
and $ESW_{A^\prime}$ such that the corresponding endomorphisms are
isospectral 

A variant of these deformations are the so called
$\sigma_A^{(a,b)}$-deformations
defined as follows.

Consider an 
$ESW_A=J_{\bold z}=\bold A\oplus\bold A^\perp$ 
such that the endomorphisms act on $\bR^n$.
For a pair $(a,b)$ of natural numbers the endomorphism
space $ESW^{(a,b)}_A=J_{\bold z}^{(a,b)}$ 
is defined by a new representation, $B^{(a,b)}=J_B^{(a,b)}$,
of the endomorphisms $B\in ESW_A$
on the new X-space 
$\bold v=\bR^n \times \dots
\times \bR^n$ 
(the 
Cartesian product is taken 
$(a+b)$-times) 
such that the endomorphisms $A^{(a,b)}$ and $F^{(a,b)}$,
where $F\in A^\perp$, are defined by

$$
\gathered
A^{(a,b)}(X)=
(A(X_1),\dots ,A(X_a),-A(X_{a+1}),\dots ,-A(X_{a+b})),\\
F^{(a,b)}(X)=
(F(X_1),\dots ,F(X_{a+b})).
\endgathered
\tag 2.13
$$

If $\sigma^{(a,b)}$ is the involutive orthogonal transformation
defined on $\bold v$ by
$$
\sigma^{(a,b)}(X)=(X_1,\dots ,X_a,-X_{a+1},\dots ,-X_{a+b}),
\tag 2.14
$$
then the 
$\sigma^{(a,b)}_A$-deforma\-tion sends
$ESW_A^{(a,b)}$ and
$ESW_A^{(a+b,0)}$ to each other.

In \cite{Sz5}
also an other type of deformations,
called $\sigma$-deformations, was introduced.
It is defined for general
endomorphism spaces 
such that the  deformation $\sigma\circ B$ is performed 
on all elements of the
endomorphism space. (Also in this case the $\sigma$ is an involutive
orthogonal transformations commuting with all the elements of the
endomorphism space.)

Though they seem to be completely different deformations,
Reduction Theorem 4.1 in \cite{Sz5} asserts 
that, on $ESW_A$'s,
$\sigma_A$-deformations are equivalent to  
$\sigma$-deformations. In this spectral investigations we prefer
the $\sigma_A$ deformations to the $\sigma$ deformations
of an $ESW_A$ because
of the simplicity offered by considering only the deformation of a
single endomorphism.

The 2-step nilpotent Lie algebras (resp. Lie groups) corresponding
to $ESW_A^{(a,b)}\allowmathbreak
=J_{\bold z}^{(a,b)}$ is denoted by 
$\frak n^{(a,b)}_J$
(resp. $N^{(a,b)}_J$).

This notation is consistent with the notation of 
$\frak h^{(a,b)}_3$.
In this case the space
$\bold z=\bR^3$ is 
identified with the space of the imaginary 
quaternions and 
the skew endomorphisms 
$J_Z=L_Z$ 
acting on $\bR^4
=\bold H$ are defined 
by left products with $Z$. Notice that in 
this case $J_{\bold z}\simeq so(3)
\subset so(4)$ hold and this endomorphism space is 
closed with respect to the Lie bracket. 

In case of Cayley numbers the Z-space
$\bold z=\bR^7$ is identified with the space of imaginary
Cayley numbers and the endomorphism space 
$J_{\bold z}=R_{\bold z}$ is defined by the right product described
below formula (2.6) (the left products result equivalent endomorphism 
spaces). Notice that this endomorphism space 
is not closed with respect to the Lie bracket.
The corresponding 
Lie algebra is denoted by
$\frak h^{(a,b)}_7$.

$\sigma$-deformations provide pairs of endomorphism spaces such that
the metrics on the corresponding groups have different local geometries
in general. 
Non-isometry Theorem 2.1 in \cite
{Sz5} asserts that
{\it 
for endomorpphism spaces $J_{\bold z}$ which are either 
non-Abelian Lie algebras 
or, more generally, they contain
non-Abelian Lie subalgebras, the metric on 
$N^{(a,b)}_J$ is locally non-isometric to the metric on
$N^{(a',b')}_J$ unless $(a,b)=(a',b')$ up to an order.
Yet the Ball$\times$Torus-type domains are both Dirichlet and Neumann
isospectral on these locally different spaces.}

The key idea of this theorem's proof is that
$\sigma^{(a,b)}$-deformations impose changing on 
the algebraic structure of the
endomorphism spaces and that is why they can not be conjugate. 

This general
theorem proves the non-isometry with respect to the groups  
$H^{(a,b)}_3$, however, it does not prove it with respect to the 
$H^{(a,b)}_7$'s, or, for
the other Cliffordian endomorphism spaces. 
Fortunately enough, the non-isometry statement in the latter case
is well known (described at (2.6)) and can be established exactly 
for those
Heisenberg-type groups, $H^{(a,b)}_l$,
where $l=3mod(4)$.

The non-isometry proofs on the solvable extensions are traced back
to the nilpotent subgroups and on the sphere-type domains 
they are traced back to the ambient space.
I. e., the question of non-isometry
is always traced back to the question of non-conjugacy of the 
corresponding endomorphism spaces and, therefore, to the above theorem.

\head
\titlebf 3. Isospectrality theorems on sphere-type manifolds
\endhead

In \cite{Sz5}, the isospectrality theorems
are completely established 
on ball- and ball$\times$to\-rus-type manifolds, however, the 
proofs are only outlined 
on the boundary, i. e., on sphere- and 
sphere$\times$torus-type manifolds. Even
these sketchy details
concentrate mostly on the striking
examples. In this chapter the isospectrality 
theorems are completely established also on these boundary manifolds.

There are three sections ahead. In the first two sections
the nilpotent case is considered where, after establishing
an explicit formula for the Laplacian on the boundary manifolds,
the isospectrality theorems 
are accomplished by 
constructing intertwining operators. In the third section these
considerations are settled on the
solvable extensions.

\head\leftline{Normal vector field and Laplacian on 
the boundary manifolds}
\endhead

We start by a brief description of the ball$\times$torus- and
ball-type domains in the 
nilpotent case.

(1) Let $\Gamma$ be a full lattice on the Z-space spanned by a basis
$\{e_1,\dots ,e_l\}$. 
For an $l$-tuple $\alpha =(\alpha_1,
\dots ,\alpha_l)$ of integers the corresponding lattice point is 
$Z_\alpha =\alpha_1e_1+
\dots 
+\alpha_le_l$. Since $\Gamma$ is a discrete subgroup, 
one can consider the
factor manifold $\Gamma \backslash N$ with the factor metric. This
factor manifold is a principal fibre bundle with the base space 
$\bold v$
and with the fibre $T_X$ at a point $X\in \bold v$. Each fibre $T_X$ is
naturally identified with 
the torus $T=\Gamma \backslash 
\bold z$. The projection 
$\pi : \Gamma \backslash N \to \bold v$ defined by
$\pi : T_X \to X$ 
projects the inner product from the horizontal subspace
(defined by the orthogonal complement 
of the
fibres) to the Euclidean inner product $\<\, ,\,\>$ on the X-space.

Consider
also a Euclidean ball $B_\delta$ 
of radius $\delta$ around the origin of the
X-space and restrict the fibre bundle onto $B_\delta$. Then the fibre
 bundle
$(B_\delta ,T)$ has the boundary $(S_\delta ,T)$, 
which is also a principal
fibre bundle over the sphere $S_\delta$.

Prior to paper \cite{Sz5} only these manifolds were involved
to constructions of isospectral metrics with different local geometries.

(2)In these papers we consider also such domains around the 
origin which are homeomorphic
to a $(k+l)$-dimensional ball and their 
smooth boundaries can be described as 
level surfaces by equations of the form
$f(|X|,Z)=0$. The boundaries of these domains are homeomorphic to
the sphere $S^{k+l-1}$ such that the 
boundary points form
a Euclidean sphere of radius $\delta (Z)$, for any fixed $Z$. 
I. e., the boundary can be
described by the equation $|X|^2-\delta^2 (Z)=0$. We
call these cases {\it Ball-cases} resp. {\it Sphere-cases}.

In this section we provide explicit formulas for  
the normal
vector field and for the Laplacian on the
sphere-type manifolds only. However, these formulas 
establish the corresponding formulas also in  
the 
sphere$\times$torus-cases, such that one substitutes 
the constant radius $R$ 
for
the function $\delta (Z)$ in order to have the formula also 
on the latter manifolds.

First 
the normal vector $\bmu$ is computed.
From the equation
$$
\nabla f=grad f=\sum_i\bold X_i(f)\bold X_i+
\sum_\al \bold Z_\al (f)\bold Z_\al
$$
we get (by using the special function 
$f(|X|,Z)=|X|^2-\delta^2 (Z)$) that this unit normal vector at a
point $(X,Z)$ is
$$
\bmu =(4|X|^2+\frac 1 {4} |J_{\nabla\delta^2}(X)|^2+
|\nabla\delta^2|^2)^
{-\frac 1 {2}}(2X-\frac 1 {2}\bold J_{\nabla\delta^2}(X)-
\nabla\delta^2),
\tag 3.1
$$
where $\bmu$ is considered as an element of the Lie algebra. 
Notice that in
the sphere$\times$torus case 
the $\bmu$ has the simple form $\bmu =X/|X|$.

By (1.2), this normal vector
can be written also in the following regular vector form:
$$
\bold \bmu=C\big(2|X|
E_0 - \sum_{\al =1}^l(\pa_\al \delta^2)
(\frac 1 {2}|X|E_\al +\sum_{\beta=1}^l(1+\frac 1 {4}\<J_\al (X),
J_\beta (X)\>)e_\beta )\big ),
\tag 3.2
$$
where $\{e_1,\dots ,e_l\}$ is an orthonormal basis of $\bold z$ and
$E_\al =J_\al (E_0);E_0=X/|X|$, furthermore,
$$
C=(4|X|^2+\frac 1 {4} |J_{\nabla 
\delta^2}(X)|^2 + |\nabla\delta^2|^2)^{-\frac 1 {2}}.
\tag 3.3
$$
In the following we make it always clear that which representation
of a particular vector is considered.

Over a fixed point $Z$, the X-crossection with the boundary $\pa D$
is the sphere $S_X(Z)$  with radius $\delta (Z)$.
The corresponding Z-crossection over a fixed point $X$ is denoted
by $S_Z(X)$. Notice that this latter manifolds are only homeomorphic
to Euclidean spheres in general and they are
Euclidean
spheres for all point $X$ 
if and only if the function $\delta$ depends only on $|Z|$.
The Euclidean(!) normal vector $\bmu_Z$ to $S_Z(X)$ is
$$
\bmu_Z=(\sum_\al (\pa_\al \delta)^2)^{-1/2}\sum_\beta \pa_\beta (\delta)
e_\beta =\sum \bmu_{Z\beta} ,
\tag 3.4
$$
which is different from the orthogonal projection of $\bmu$ onto the
Z-space. 

Let $\tilde{\na}$ (resp. $\tilde\Delta$) be the covariant derivative
(resp. the Laplace operator) on the boundary $\pa D$.
The second fundamental form and the Minkowski curvature are denoted by
$M(V,W)$ and $\bold M$. Then the formula
$$
\na^2f(V,V)=V\cdot V(f)-\na_VV\cdot (f)=\tilde{\na}^2f(V,V)+M(V,V)
f^\prime 
\tag 3.5
$$
($f^\prime :=\bmu\cdot (f)$) holds for any function $f$ defined on
the ambient space and for any vector field $V$ tangent to $\pa D$.
Thus
$$
\tilde{\Delta}f=\Delta f-f^{\prime\prime}-\bold Mf^{\prime} .
\tag 3.6
$$
Choose such functions $f$ around the boundary 
$\pa D$ which are constant with
respect to the normal direction (i. e. $f^\prime =f^{\prime\prime}=0$).
Then check the formula 
$$
\gathered
\tilde \Delta =\Delta_{S_X(Z)}+
\Delta_{S_Z(X)}+ 
\sum_{\al =1}^l(\pa_\al -\bmu_{Z\al})
D_\al\bullet +
\\
\frac 1 {4}\sum_{\al\beta =1}^l\<J_\al (X),J_\beta (X)\> 
(\pa_\al -\bmu_{Z\al})(\pa_\beta -\bmu_{Z\beta}).
\endgathered
\tag 3.7
$$

This
formula is simpler on the level surfaces described by
equations of the form $|X|=\delta (|Z|)$ 
(cf. Chapter 4).
On the sphere$\times$torus-type manifolds we get the Laplacian by
performing
the simple modification $\bmu_{Z\al}=0$ in the above formula.

\head
\leftline{Isospectrality theorems on sphere-type manifolds} 
\endhead

In order to establish 
the isospectrality theorems 
on sphere-type manifolds, one should appropriately modify the technique
developed for the ambient ball-type
manifolds in \cite{Sz5}. 
The main tool of this technique is the following
 
\head\leftline{\it Harmonic analysis developed for a unit
anticommutator $J_A$}
\endhead

A brief account of this analysis is as follows.

As it is indicated, the anticommutator $A$ is a unit
anticommutator. By the normalization described at the beginning of
\S 2, each non-degenerated anticommutator 
can be rescaled to a unit anticommutator. 

First notice that the Euclidean Laplacian $\Delta_S$ (defined 
on the unit sphere S
around the origin of the X-space) and the differential operator 
$\bold D_A\bullet$ commute since the vector
field $J_A(X_u)$, where $X_u\in S$, is an infinitesimal generator of 
isometries on $S$. Therefore a 
common eigensubspace decomposition of
the $L^2$ function space exists which can be established 
as follows.

The eigenfunctions of $\Delta_S$ are the well known spherical
harmonics which are the restrictions of the homogeneous 
harmonic polynomials
of the ambient X-space onto the sphere S. The space of the
$q^{th}$-order
spherical harmonic polynomials is denoted by $\bold H^{(q)}$.
In the following we describe the eigenfunctions of the operator
$\bold D_A\bullet$.

For a fixed X-vector
$Q$, we define the complex valued function
$$
\bold \Theta_Q(X)=\<Q+\bold iJ_A(Q),X\> =\<\bold Q,X\>,
\tag 3.8
$$
where 
$\bold Q=Q+\bold iJ_A(Q)$.
Then the polynomials of the form
$$
\Phi_{(Q_i,p_i,Q^*_j,p^*_j)}(X)=\Pi_{i=1}^r\bold \Theta_{Q_i}^{p_i}(X)
\Pi_{j=1}^{r^*}\overline {\bold \Theta}_{Q^*_j}^{p_j^*}(X)
\tag 3.9
$$
are eigenfunctions of the operator $\bold D_A\bullet$ with the
eigenvalue $(2s-q)\bold i$, where
$q=\sum_i p_i+\sum_jp_j^*=s+(q-s)$. 

Notice, that the functions of the pure form
$$
\Phi_{(Q_i,p_i)}(X)=\Pi_{i=1}^r\bold \Theta_{Q_i}^{p_i}(X)
\quad ;\quad \overline{ 
\Phi}_{(Q_i,p_i)}(X)=\Pi_{i=1}^r\overline {\bold \Theta}_{Q_i}^{p_i}(X)
\tag 3.10
$$
are harmonic with respect to the Euclidean Laplacian $\Delta_X$ 
on the $X$-space. In fact, on the Euclidean K\" ahler manifold
$\{ \bold v,\< ,\> ,A\}$ 
these functions correspond to the holomorphic
resp. anti-holomorphic polynomials. One can directly 
check this property also
by $\<\bold Q_i,\bold Q_j\>=
\<\overline{\bold Q_i},\overline{\bold Q_j}\>=0$.

However, the polynomials of the mixed form are not
harmonic, since 
$$ 
\Delta_X\Theta_Q(X)\overline Q_{Q^*}(X)=2\<Q,Q^*\>+2\bold i
\< J_A(Q),Q^*\>.
\tag 3.11
$$ 
Let us note that
the whole space $\bold H^{(q)}$ of $q^{th}$ order eigenfunctions is 
not spanned
by the above polynomials of pure form. 
The "missing" functions can be furnished by the orthogonal projection
of the $q^{th}$ order mixed polynomials $\Phi_{(Q_i,p_i,Q^*_i,p_i^*)}$
 onto the function space 
$\bold H^{(q)}$. The
range of this projection is denoted by $\bold H^{(s,q-s)}$. 
Since the operators $\Delta_S$ and
$D_A\bullet$ commute, 
the space $\bold H^{(q)}$ is invariant under
the action of $D_A\bullet$ and the subspace 
$\bold H^{(s,q-s)}\subset \bold H^{(q)}$
is an eigensubspace of this operator with the eigenvalue 
$(2s-q)\bold i$.
Thus the decomposition 
$\bold H^{(q)}=\oplus_{s=0}^q\bold H^{(s,q-s)}$ is an 
orthogonal direct sum corresponding to the common 
eigensubspace decomposition
of the two commuting differential operators $\Delta_S$ and 
$D_A\bullet$.

A more accurate description of the above mentioned projections
can be given by
the 
kernel functions $H_{(q)}(Q_u,Q^*_u)$ introduced for the
subspaces $\bold H^{(q)}$ by
$$
H_{(q)}(Q_u,Q^*_u)=\sum_{j=1}^{N_q}\eta_j^{(q)}(Q_u)\eta_j^{(q)}(Q^*_u),
\tag 3.12
$$
where $\{\eta_1^{(q)},\dots ,\eta_{N_q}^{(q)}\}$ is an orthonormal basis
on the subspace $\bold H^{(q)}$. In \cite {Be} it is proved (cf. 
Lemma 6.94) that 
the eigenfunction $H_{(q)}(Q_u,.)$
is radial 
for any fixed $Q_u$. 
(I. e., it has the form
$C_q\<Q_u,.\>^q+\dots +C_1\<Q_u,.\> +C_0$ with
$H_{(q)}(Q_u,Q_u)=1$.)
Furthermore, for any function $\Psi\in L^2(S)$ we have 
$$
\Psi(Q_u)=\sum_{q=0}^\infty\int H_{(q)}(Q_u,Q^*_u)\Psi (Q^*_u) dQ^*_u=
\sum h_{(q)}(\Psi )_{/Q_u},
\tag 3.13 
$$
which is called 
the {\it spherical
decomposition of $\Psi$ by the spherical harmonics}. The operators
$
h_{(q)} :L^2(S)\to\bold H^{(q)}
$
project the $L^2$ function space to the corresponding eigensubspace of 
the Laplacian.

For an $r^{th}$ order polynomial 
$\Pi_i \Theta_{Q_i}^{p_i}\overline{\Theta}^{r_i-p_i}_{Q_i}$,
$\sum r_i=r$, 
the projection $h_{(r)}$ can be computed also by the formula
$$
h_{(r)} 
(\Pi_i \Theta_{Q_i}^{p_i}\overline{\Theta}^{r_i-p_i}_{Q_i})
=\sum_s B_s \<X,X\>^s\Delta^s_X\Pi_i 
\Theta_{Q_i}^{p_i}\overline{\Theta}^{r_i-p_i}_{Q_i},
\tag 3.14
$$
where $B_0=1$ and the other coefficients can be determined
by the recursive formula
$$
2s(2(s+r)-1)B_s+B_{s-1}=0.
$$
These formulas are established by the fact that the function
on the right side of (3.14) is a homogeneous harmonic polynomial
exactly for these 
coefficients.

One of the most important properties of these operators is that 
they commute
with the differential operators $D_\al\bullet$. This statement
immediately follows from the fact that the vector fields 
$J_\al (X)$ are 
infinitesimal generators of one parametric families of 
isometries on the Euclidean 
sphere $S$ and the projections $h_{(q)}$ are invariant
with respect to these isometries. One can imply this commutativity
also by (3.14) and by the commutativity of the operators $\Delta_X$
and $D_\al\bullet$.

By substituting $Q_u=1/2(\bold Q_u+\overline {\bold Q}_u)$,
where $\bold Q=Q+\bold{i}J_A(Q)$, into the above
expression of the radial kernel $H_{(q)}(Q_u,.)$, we get
\proclaim {Proposition 3.1 \cite {Sz5}} 
The $r^{th}$ order polynomial space $\bold
P^{(r)}$ is the direct
sum 
of the subspaces 
$\bold P^{(p,r-p)}$ spanned by the polynomials of the form
$\Theta_{Q}^p\overline {\Theta}_{Q}^{r-p}$, where $Q\in\bold v$
and $0\leq p\leq r$. The space 
$\bold P^{(p,r-p)}$ consists of the $r^{th}$
order eigen polynomials 
of the differential
operator $D_A\bullet$
with eigenvalue $(r-p)\bold i$.

The projection 
$h_{(r)}$
establishes a one to one and onto map between the polynomial spaces
$\bold P_S^{(p,r-p)}$ and $\bold H^{(p,r-p)}$. 
We denote this restricted
map by 
$h_{/(r)}: \bold P_S^{(r)}\to\bold H^{(r)}$.
The direct sum 
$T=\oplus_rh_{/(r)}$ of these maps
defines an invertible operator on the whole function space 
$L_S^2$. This operator $T$ is commuting with the differential operators
$D_\al\bullet$.
\endproclaim 

Above,
$\bold P_S^{(p,r-p)}$ denotes the space of 
the corresponding restricted functions onto the sphere $S$. 
In the following we represent the functions $\psi\in\bold H^{(r)}$
in the form $\psi=h_{/(r)}(\psi^*)$, where $\psi^*\in\bold P^{(r)}_S$.

On 
the whole ambient space $\bold n$, the function
space is spanned by the functions of the form
$$
F(X,Z)=\varphi (|X|,Z)h_{/(r)}(\Theta^p_{Q_u}
\overline {\Theta}_{Q_u}^{r-p})(X_u),
\tag 3.15
$$
where $X_u=X/|X|$ is a unit vector. The decomposition with respect to
these 
functions is called {\it spherical decomposition}. 
If we cancel $h_{/(r)}$
in (3.14), the corresponding decomposition 
with respect to these functions
is called {\it polynomial decomposition}. The operator $T$ transforms
a polynomial decomposition to a spherical decomposition.

\head\leftline
{\it Constructing the intertwining operator}
\endhead

Let $J_{\bold z}=J_{\bold A}\oplus J_{\bold A^\perp}$ and
$J_{\bold z^\prime}=J_{\bold A^\prime}\oplus J_{\bold A^{\prime\perp}}$
be endomorphisms spaces with the unit anticommutators 
$J_A$ and $J_{A^\prime}$ such that the endomorphisms 
are acting on the same space and, even more, 
$J_{\bold A^\perp}=J_{\bold A^{\prime
\perp}}$ holds. In \cite{Sz5} we proved that the
ball- and ball$\times$torus-type domains with
the same radius-function 
are both Dirichlet and Neumann isospectral on the metric groups 
$N=N_{J_{\bold z}}$ and $N^\prime =N_{J_{\bold z^\prime}}$
constructed by these $ESW_A$'s. This isospectrality 
theorem is established
for $\sigma_A$-deformation 
where $A$ is not necessarily a unit anticommutator.

For the proof of this theorem we constructed the intertwining operator
$$
\gathered
\kappa :L^2(N,\bold C)\to 
L^2(N^\prime ,\bold C) \\
\kappa :
F(X,Z)\to
F^\prime (X,Z)=\varphi (|X|,Z)h^\prime_{/(r)}
(\Theta ^{\prime p}_{Q_u}\overline{\Theta ^\prime}^{r-p}_{Q_u})
(X_u),
\endgathered
\tag 3.16
$$
where $F$ is introduced in (3.15) and the functions $\Theta^\prime$ are
defined on $N^\prime$ by means of 
$J_{A^\prime}$.

This map maps the function space $\bold H^{(r)}$ onto
$\bold H^{\prime (r)}$ and it is defined by means of the map
$$
\kappa^*:\bold P^{(r)}\to\bold P^{\prime (r)}\quad ,\quad
\kappa^*:
\Theta^p_Q\overline{\Theta}_Q^{r-p}\to
\Theta ^{\prime p}_Q\overline{\Theta} _Q^{\prime r-p}
\tag 3.17
$$
and of the projection $h_{/(r)}$. Unlike $\kappa$, 
the map $\kappa^*$ can be easily handled.
If $\{E_1,\dots ,E_K\}$ is a basis on the X-space then
$\kappa^*(\Pi\Theta_{E_{i_r}}\Pi\overline{\Theta}_{E_{j_r}})=
\Pi\Theta^\prime_{E_{i_r}}\Pi\overline{\Theta}^\prime_{E_{j_r}}$.
In general, for arbitrary vectors $Q_m$, we get
$$
\kappa^*(\Pi\Theta_{Q_{i_r}}\Pi\overline{\Theta}_{Q_{j_r}})=
\Pi\Theta^\prime_{Q_{i_r}}\Pi\overline{\Theta}^\prime_{Q_{j_r}}.
\tag 3.18
$$

The operators $\kappa$ and $\kappa^*$ are connected by the equation
$\kappa=T^\prime\circ\kappa^*\circ T^{-1}$.
In \cite{Sz5}, the intertwining property of the map 
$\kappa$ is proved. 

The intertwining operator $\pa\kappa$
on the Sphere-type boundary is construct\-ed by 
an appropriate restriction
of $\kappa$ onto the boundary.
Also in this case, first we suppose that the 
$J_A$ is a unit anticommutator and in the end we make the necessary
modifications in order to establish the intertwining 
for general $\sigma_A$-deformations.

The function space $L^2(\pa B)$ is spanned
by the functions of the form 
$$
F(\pa X,\pa Z)=\varphi (|\pa X|,\pa Z)h_{/(r)}(\Theta^p_{Q_u}
\overline {\Theta}_{Q_u}^{r-p})(X_u),
\tag 3.19
$$
where $(\pa X,\pa Z)\in \pa B$. Then the operators 
$\pa\kappa ,\pa\kappa^* L^2(\pa B)\to L^2(\pa B)$ 
are defined again by the formulas (3.16)-(3.17), what are used for
defining $\kappa$ and $\kappa^*$ on the ambient space. However, 
in this case,
the function $\varphi$ depends on variables $(|\pa X|,\pa Z)$. 
Then the intertwining property can be proved by the
following steps. 

Consider the explicit expression
(1.6) of the Laplacian $\tilde\Delta$ on $\pa B$.
Since $\pa \kappa :\bold H^{(q)}\to \bold H^{\prime (q)}
=\bold H^{(q)}$, 
the terms $\Delta_{S_X(Z)} + \Delta_{S_Z(X)}$ and 
$\Delta^\prime_{S_X(Z)} +\Delta^\prime_{S_Z(X)}$ are 
clearly intertwined by this map. 

To the next step choose an orthonormal basis 
$\{e_0,e_1,\dots ,e_{l-1}\}$ on the Z-space
such that $e_0=A$. The Greek characters are used for the indices
$\{0,1,\dots ,l-1\}$ and the Latin characters are used for the indices
$\{1,\dots ,l-1\}$. Then
$D_c \bullet =D^\prime_c \bullet$
furthermore
$$
D_0\bullet \bold \Theta_Q(X)=\bold i\bold \Theta_Q(X)\quad ,\quad 
D_0\bullet \overline {\Theta}_Q(X)=-\bold i\overline {\Theta}_Q(X),
\tag 3.20
$$
and
$$
D_c\bullet
\Theta_{Q}=\<\bold Q,J_c(X)\>=
-\overline{\Theta}_{J_c(Q)}\quad ,\quad
D_c\bullet\overline{\Theta}_{Q}=-
\Theta_{J_c(Q)}.
\tag 3.21
$$
(Let us mention that the switching of the conjugation in (3.20)
is due to the equation $J_0J_c=-J_cJ_0$).
Therefore 
$\pa\kappa^*D_\al\bullet (\psi^*)
=D'_\al\bullet
\pa\kappa^*(\psi^*)$
for any $\psi^*\in\bold P^{(r)}_S$. Since 
$\pa\kappa =T'\circ\pa\kappa^*\circ T^{-1}$, we also have
$\pa\kappa D_\al\bullet (\psi)
=D'_\al\bullet
\pa\kappa(\psi)$ and thus
the terms $\sum (\pa_\al -\bmu_{Z\al}) D_\al \bullet$ and 
$\sum (\pa_\al -\bmu_{Z\al}) D'_\al \bullet$
in (1.6)
are intertwined by the map $\pa\kappa$.

Only the term $(1/4)\sum \<J_\al (X),J_\beta (X)\>
(\pa_\al-\bmu_{Z\al})(\pa_\beta -\bmu_{Z\beta})$ should be 
considered yet.

First notice that on the Heisenberg-type 
groups $H^{(a,b)}_l$, this
operator is nothing but $(1/4)|X|^2\Delta_{S_Z(X)}$. Therefore it is 
intertwined by the $\pa\kappa$ and the proof is completely 
established on this
rather wide range of manifolds defined by Clifordian
endomorphism spaces $J_l^{(a,b)}$.

Now consider this operator in general cases.
Since $J_0\circ 
J_c$ is a skew symmetric endomorphism, therefore
$\<J_0(X),J_c(X)\>=0$. Thus this 
operator is the same one on the considered spaces 
$N$ and $N^\prime$.
Yet, for the sake of completeness, we should prove that the 
$\pa\kappa$ maps
a function of the form 
$$
\<J_c(X),
J_d(X)\>h_{(r)}(\Theta^p\overline
{\Theta}^{r-p}):=J_{cd}(X)h_{(r)}(\Theta^p\overline{\Theta}^{r-p})
\tag 3.22
$$ 
to a function of the very same form on $N^\prime$.
In \cite{Sz5} 
this problem is settled 
(cf. formulas (4.13)-(4.17)) by proving, first, that
in the spherical 
decomposition of the function 
$J_{cd}\Theta_Q^p\overline{\Theta}^{r-p}$
the component-spherical-harmonics 
are linear combinations of the functions of the form
$$
\Theta_P\overline{\Theta}_R
\Theta^{p-s}_{Q_u}\overline{\Theta}^{r-p-s}_{Q_u}\quad ,\quad
\Theta^{p-v}_{Q_u}\overline{\Theta}^{r-p-v}_{Q_u}
$$
such that the combinational coefficients depends only on
the constants 
$$
r\, ,\,p\, ,\,s\, ,\,v\, ,\,Tr(J_c\circ
J_d)\, ,\,\<J_c(Q_u),J_d(Q_u)\> .
$$
Then the same statement is established for the preimages (with respect
to the map $T$) of these functions. 
Since these terms do not depend on $J_A$,
the function
$J_{cd}h_{(r)}\Theta_Q^p\overline{\Theta}^{r-p}_Q$ is
intertwined with the corresponding function
$J_{cd}h_{(r)}\Theta^{\prime p}_Q\overline{\Theta}^{\prime r-p}_Q$
by the map $\kappa$. 

This argument settles the proof of intertwining also for $\pa\kappa$,
yet we would like to give a more simple and comprehensive proof for 
this part here.

The radial spherical harmonics span the eigensubspaces $\bold H^{(r)}$
(cf. (3.13)). Therefore, it is enough to consider the real functions 
of the form $K_{(r)}(X)=J_{cd}(X)h_{/(r)} \< Q,X\>^r$.
In this case the function
$h_{/(r)} \< Q,X\>^r$ is nothing but a constant multiple of
$H_{(r)}(Q,X)$. By (3.14), the spherical decomposition of $K_{(r)}$ is
$K_{(r)}=\sum C_p\< X,X\>^pL_p(X)$, where the functions $L_p(X)$ are
spherical harmonics built up by the functions of the form

$$
\gathered
\<Q,X\>^p\, ,\,
J_{cd}(X)\<Q,X\>^q\, ,\\
 \<(J_cJ_d+J_dJ_c)(Q),X\>^s
\<Q,X\>^{v}=\<Q_{cd},X\>^s  
\<Q,X\>^{v},
\endgathered
\tag 3.23
$$  
such that the combinational coefficients depends only on
the constants $r,p,s,v,\allowmathbreak TrJ_c\circ
J_d$ and on $\<J_c(Q),J_d(Q)\>$. More precisely,
the function $L_p$ is of the form $L_p=h_{(p)}U_p$, where $U_p$ is one
of the functions from the set (3.23).

The function $J_{cd}(X)$ 
can be written in the form
$$
\gathered
J_{cd}(X)=
\sum_{i=1}^K \frac 1 {4} \{ (
\< X,\bold Q_{ci}\>\< X,\overline {\bold Q}_{di}\> +\< X,\overline
{\bold Q}_{ci}\>\< X,\bold Q_{di}\> ) + \\ (
\< X,\bold Q_{ci}\>\< X,\bold Q_{di}\> +\< X,\overline
{\bold Q}_{ci}\>\< X,\overline{\bold Q}_{di}\> )
\} = J^{(1)}_{cd}(X)+J^{(2)}_{cd}(X),
\endgathered
\tag 3.24
$$
where $E_1,\dots ,E_K$ is an orthonormal basis on the X-space 
($K=k(a+b)$) and
$\bold Q_{ei}=\bold J_e(\bold E_i)$. 
The proof of this formula immediately follows from 
$$
J_{cd}(X)=\sum_i\<J_c(X),E_i\>\<J_d(X),
E_i\>\quad ,\quad E_i
=\frac 1 {2} (\bold E_i+\overline{\bold E}_i).
\tag 3.25
$$

Notice that the second function of (3.24) is vanishing.
This statement immediately follows from equations
$\<\bold Q_1,\bold Q_2\>=
\<\overline{\bold Q}_1,\overline{\bold Q}_2\>=0$.

By the substitution $Q=(1/2)(\bold Q+\overline{\bold Q})$
we get that the spherical harmonics $L_p$ 
are linear combinations of the functions of the form
$$
\Theta^{p-s}_{Q}\overline{\Theta}^{r-p-s}_{Q}\, ,\,
J_{cd}
\Theta^{p-s}_{Q}\overline{\Theta}^{r-p-s}_{Q}\, ,\,
(\Theta_{Q_{cd}}+\overline{\Theta}_{Q_{cd}})
\Theta^{p-v}_{Q}\overline{\Theta}^{r-p-v}_{Q}
\tag 3.26
$$
such that the combinational coefficients depends only on
the constants $r,p,s,v,\allowmathbreak TrJ_c\circ
J_d$ and $\<J_c(Q),J_d(Q)\>$.

The considered problem can be settled by representing $J_{cd}$ in the
form (3.24) in formula (3.22).
In fact, the coefficients discussed above
are the same on both spaces (since they
do not depend on the unit anticommutator $J_0$). By
(3.14) we get that
the preimages 
(with respect to the map $T$) 
of spherical-harmonics $L_p$
are
combinations of appropriate 
functions which have the same form (3.26) on the other manifold and
the coefficients do not depend on
$J_0$. This completely proves
that $\pa\kappa$ maps a function $K_{(r)}$ 
to an appropriate function
desired in this problem.

The above constructions and proofs can be easily extended to the 
cases, when $J_A$ is just a non-degenerated anticommutator and
the anticommutators $A$ and $A^\prime$ are $\sigma_A$-related.
 
In this case, first,
the unit anticommutator $A_0$ should be introduced by 
an appropriate rescaling of $A$.
This $A_0$ may not be in the endomorphism space, however, 
it is commuting
with $A$ and it is anticommuting with the elements of 
$J_{\bold A^\perp}$.
The operator $\pa\kappa$ should be established by means of $A_0$.

Since
$\sigma$-deformations do not change the maximal eigensubspaces
of the operators involved, 
the operators $D_{\alpha}\bullet$ and $D_{\alpha}^{\prime}\bullet$
are intertwined by$\pa\kappa$. By the very same reason also 
the operators
$\<J_A(X),J_A(X)\>(\pa_0 -\bmu_{Z0})^2$ and
$\<J_{A^\prime}(X),J_{A^\prime}(X)\>(\pa_0 -\bmu_{Z0})^2$ are 
intertwined by the $\pa\kappa$. 
Since these are the only terms in the Laplacian
which depend on the eigenvalues of the endomorphism $J_A$, the
intertwining property is established in the considered case. 
Thus we have.

\proclaim{Main Theorem 3.2} Let 
$J_{\bold z}=J_{\bold A}\oplus J_{\bold A^\perp}$ and
$J_{\bold z^\prime}=J_{\bold A^\prime}\oplus J_{\bold A^{\prime\perp}}$
be endomorphism
spaces 
acting on the same space such that  
$J_{\bold A^\perp}=J_{\bold A^{\prime\perp}}$, furthermore, the 
anticommutators $J_A$ and $J_{A^\prime}$ are either unit endomorphisms
or they are $\sigma$-related.
Then the map 
$\pa\kappa =T^\prime\circ\pa\kappa^* T^{-1}$ intertwines the 
corresponding Laplacians on the sphere-type boundary $\pa B$
of any ball-type domain on the metric groups $N_J$ and $N_{J^\prime}$.
Therefore the corresponding metrics on
these sphere-type manifolds are isospectral.
\endproclaim
\medskip\noindent
{\bf Remark 3.3}
The map $\pa\kappa$ establishes the
isospectrality theorem also on the sphere$\times$torus-type
boundaries of the ball$\times$torus-domains in the considered cases,
offering a 
completely new proof for the theorem.
 
\head\leftline{Intertwining operators 
on the solvable extensions}
\endhead

The above isospectrality theorem extends to the 
solvable extensions of nilpotent groups. In this case one should
consider the following domains.

{\it The Solvable Ball $\times$ Torus Cases:} A group $SN$ can be 
considered as a
principal fibre bundle (vector bundle) 
over the $(X,t)$-space, fibrated by
the Z-spaces. Let 
$\Gamma$ be again a full lattice on the Z-space, and $D_{R(t)}$ be
a domain on the $(X,t)$-space such that it
is diffeomorphic to a $(k+1)$-dimensional
ball whose smooth boundary can be described by an equation of the form
$|X|=R(t)$. Then consider the 
torus bundle 
$(\Gamma /\bold z,D_{R(t)})$ 
over $D_{R(t)}$.
The normal vector $\bold \bmu$ at a boundary point 
$(X,Z,t)$ is of the form
$\bold \bmu=A(t)X+B(t)\pa_t$, where $A(t)$ and $B(t)$ are determined by 
$R(t)$.

{\it The Solvable Ball Case:} In this case we consider a domain 
$D$ diffeomorphic
to a $(k+l+1)$-dimensional ball whose smooth boundary (diffeomorphic to
$S^{k+l}$) can be described  as a levelset in the form 
$|X|=\delta (Z,t)$.

The normal vector $\bmu$ at a boundary point 
$(X,Z,t)$ can be similarly computed
as in the nilpotent case. Then, by (1.10) and (3.1),
we get:
$$
\bmu =F(|X|,Z,t,c)t^{\frac 1{2}}
(2X-\frac 1 {2}\bold J_{grad_Z\delta^2}(X))-t(
grad_Z\delta^2 +c\pa_t (\delta^2)\bold T),
\tag 3.27
$$
where $\bmu$ is considered as an element of the Lie algebra and
$$
F=
(4t|X|^2+\frac 1 {4}t|J_{grad_Z\delta^2}(X)|^2+t^2(
|grad_Z\delta^2|^2 +c^2(\pa_t\delta^2)^2))^
{-\frac 1 {2}}.
\tag 3.28
$$
Therefore this vector can be written in the following regular vector
form
$$
\gathered
\bmu = F_0(|X|,Z,t,c)E_0
+C(Z,t,c)\pa_t + \\
\sum_{i=1}^l(F_i(|X|,Z,t,c)E_i+L_i(|X|,Z,t,c)e_i),
\endgathered
\tag 3.29
$$
where the functions $F_\al\, ,\, L_i$ 
and $C$ are determined by 
$\delta (Z,t)$ and $c$.

The Laplacian can be established by formulas (1.12),(3.5) and (3.6).
Then we get
$$
\gathered
\tilde
\Delta =t\Delta_{S_X}+t^{\frac 1 {2}}\Delta_{S_Z}+ \\
\frac 1 {4}t\sum_{\al ;\beta =1}^l
\<J_\al (X),J_\beta (X)\>(\pa_\al -\bmu_{Z\al})
(\pa_\beta -\bmu_{Z\beta}) \\
+t\sum_{\al =1}^l(\pa_\al -\bmu_{Z\al}) D_\al\bullet+
c^2t^2(\pa_t-\bmu_t)^2 +
c^2(1-\frac k {2}-l)t(\pa_t -\bmu_t).
\endgathered
\tag 3.30
$$

By repeating the very same arguments used in the nilpotent cases
(only the function $\varphi$ in (3.19) should be of the
form $\varphi (|\pa X|,\pa Z,t)$), we get

\proclaim {Main Theorem 3.4} Let $J_{\bold z}$ and $J_{\bold z^\prime}$
be endomorphism spaces described in Main Theorem 3.2. Then the 
corresponding metrics on the 
sphere-type surfaces having the same radius-function $\delta$ 
are isospectral on the solvable groups
$SN_J$ and $SN_{J^\prime}$.
\endproclaim

\head 
\titlebf 4. Extension and  
non-isometry theorems on sphere-type manifolds 
\endhead

The non-isometry theorems are established by an independent statement 
asserting that an isometry between 
two sphere-type manifolds extends to an isometry between the 
corresponding ambient manifolds. Therefore the
non-isometry  
on the sphere-type boundary manifolds can be checked by checking
the non-isometry on the ambient manifolds.
Since the non-isometry proofs on the solvable ambient manifolds are
traced back to the nilpotent cases, where the question
of non-isometry is equivalent to the non-conjugacy of the endomorphism
spaces involved, one can always check on the non-isometry simply
by checking the non-conjugacy of the corresponding endomorphism
spaces. These kind of theorems, concerning the non-conjugacy of
spectrally equivalent endomorphism spaces, are 
established in \cite{Sz5} and are reviewed in
the last part (cf. below formula) of Section 2 in this paper.
 
The proof of this {\it extension theorem} is rather complicated, due 
to the circumstances that
no general technique
has been found covering  
the diverse  
isospectrality examples constructed
in this paper. 
It requires 
different techniques depending on the sphere-type manifolds.
On the largest class of examples the scalar curvature is used
such that the extension of an isometry
from a sphere-type boundary to 
the ambient space is settled
for those manifolds where the gradient of the
scalar curvature is non-vanishing almost everywhere (this assumption is
formulated in a more precise form later). 
However, this proof does not cover the important
case of the striking examples, since the considered geodesic spheres
have constant scalar curvature. In this case
the Ricci curvature should be involved to establish the
desired extension of the isometry onto the ambient space. 
On this examples
we establish more 
non-isometry proofs,
revealing surprising spectrally undetermined objects. 
The most surprising revelation is that the spectrum of the Laplacian
acting on functions may have no information  
about the isometries.

The proofs are described in a hierarchic order. First in 
the nilpotent-  and then in
the solvable-case such sphere-type manifolds are considered
which satisfy the above mentioned condition concerning 
the scalar curvature. The striking examples are considered in
the third part. There are hierarchies also within these groups 
of considerations.
For instance, in the first big group of the proofs, first the 
Heisenberg-type nilpotent groups are considered, since
they provide a simple situation. Yet this proof clearly points into the 
direction of a general solution.

\head\leftline{\bf The nilpotent case}
\medskip\leftline{Technicalities on sphere-type manifolds}
\endhead

In order to avoid long technical computations, 
we give detailed extension and non-isometry proofs on the particular
sphere-type domains
which can be described as level sets by equations of the form
$\varphi (|X|,|Z|)=0$. (For local description we use the explicit
function of the form $|X|=\delta (|Z|)$ or an appropriate variant of
this function.) For such level sets both the X-crossections, $S_X(Z)$,
over a point Z and the Z-crossections, $S_Z(X)$, over a point X are
Euclidean spheres in the corresponding Euclidean spaces.

Without proof let us mention
that the geodesic spheres around the origin of a
Heisenberg-type group belong to this category. (We do not use this
Statement in the following considerations.  Later we independently prove
that the geodesic spheres around the origin
$(0,1)$ of the solvable extension $SH$ of a Heisenberg-type group
$H$ are level sets of the form $\varphi (|X|,|Z|,t)=0$.)

The normal vector $\bmu$ can be computed by means of formulas (1.2). 
By these formulas, regular X- and Z-vectors
tangent to the sphere-type manifold
can be expressed as Lie algebra
elements. After such a computation we get that
the perpendicular normal vector has the
form: 
$$
\bmu =
C(2X-D^\prime \bold J_Z(X))-
2CD^\prime Z=\bmu_X +\bmu_ Z,
\tag 4.1
$$
where $\bmu$ is considered as 
an element of the Lie algebra. The function
$D$ is defined by
 $D(|Z|^2)=\delta^2(|Z|)$, furthermore
$$
C=(4|X|^2+ (D^\prime )^2(|\bold J_Z(X)|^2
+4|Z|^2))^{-\frac 1 {2}}.
\tag 4.2
$$
(The prime in formula
$D^\prime$ means differentiation with respect to
the argument $\tau =|Z|^2$.)

The Weingarten map $B(\tilde U)=\nabla_{\tilde U}\bmu$ or the
second fundamental form $M(\tilde U,\tilde {U}^*)
\allowmathbreak =
g(\tilde U,B(\tilde{U}^*))$ ,
where $\tilde U$ and $\tilde U^*$ are tangent to the hyper-surface, 
can be computed by the decompositions $X=\sum x^i\bold X_i$ and 
$Z=\sum z^\al\bold Z_\al$. Then, by (1.2) and (1.3), we get:
$$
\aligned
M(\wt{X}_1,\wt{X}_2)& =C
(2\<\wt{X}_1,\wt{X}_2\>-
\sum_\beta d_\beta
\<J_\beta (X),\wt X_1\>\<J_\beta (X),\wt X_2\>); \\
M(\wt{Z}_1,
\wt{Z}_2)& =-2CD^\prime \<\wt{Z}_1,\wt{Z}_2\> ;\\
M(\wt{X},\wt{Z})& =
M(\wt{Z},\wt{X})=-\frac 1{2}\<J_{\wt{Z}}(\bmu_X+2CD^\prime X),
\wt{X}\>,
 \endaligned
\tag 4.3
$$
where $d_0={1\over2}D^\prime +|Z|^2D^{\prime\prime}$ and $d_i={1\over2}
D^\prime ,\forall\,\, 0<i\leq (l-1)$. In the first formula
an orthonormal basis 
$e_0,e_1,\dots ,e_{l-1}$ is considered 
on the Z-space such that $e_0=Z/|Z|$.

The Riemannian curvature of the considered hypersurfaces can be 
computed by the Gauss equation:
$$\aligned
&\wt R( \wt V, \wt Y) \wt W= R(\wt V, \wt Y) \wt W
-\langle R(\wt V,\wt Y)\wt W,\bmu\rangle\bmu\\
&\qquad+ M(\wt Y, \wt
W) B(\wt V)- M(\wt V, \wt W) B(\wt Y).
\endaligned
\tag 4.4
$$
In the following we compute also the Ricci curvature $\wt r (\wt U,
\wt V)$ and the scalar curvature $\wt \kappa =Tr(\wt r)$
on these
hypersurfaces. From the Gauss equation and from (1.14) we get:
$$
\aligned
\wt r(\wt Y,\wt W) & = r(\wt Y,\wt W)-\langle R(\bmu,\wt Y)
\wt W,\bmu\rangle\\
&\qquad+\big\langle\wt Y,\big((\Tr B)B-B^2\big)
(\wt W)\big\rangle;\\
\langle R(\bmu,\wt X_1),\wt X_2,\bmu\rangle
&=-{3\over4}\sum_\al\langle \wt X_1,J_\al (\bmu_X)\rangle\langle
\wt X_2,J_\al (\bmu_X)\rangle\\ &\qquad
+\frac 1{4}\langle J_{\bmu_Z}
(\wt X_1),J_{\bmu_Z}(\wt X_2)\rangle;\\
\langle R(\bmu,\wt Z)\wt X,\bmu\rangle&=\frac 1{2}\langle (J_{\bmu_Z}
J_{\wt Z}-\frac 1{2}J_{\wt Z}J_{\bmu_Z})(\wt X),\bmu_X\rangle;\\
\langle R(\bmu,\wt Z_1),\wt Z_2,\bmu\rangle&={1\over4}
\langle J_{\wt Z_1}(\bmu_X),J_{\wt Z_2}(\bmu_X)\rangle.
\endaligned
\tag 4.5
$$

In the following the scalar curvature $\wt\kappa$ 
is used to establish the 
extension theorem. In general, this 
scalar curvature is a complicated
expression depending on the functions 
$$
\aligned &
|Z|\, ,\, |J_Z(X)|\, ,\, |J^2_Z(X)|\, ,\, \sum_i |J_i(X)|^2\, ,\,
\sum_{i,j} \<J_i(X),J_j(X)\>,\\&\qquad
\sum_i \<J_i (X),J_Z(X)\>\, ,\,
\sum_i \<J_i J_Z(X),J_i J_Z(X)\>,
\endaligned
\tag 4.6
$$
where the index $0$ concerns the unit vector $Z_0$ and the indices
$i,j\,>0$ concern an orthonormal system $\{Z_1,\dots ,Z_{l-1}\}$ of
vectors perpendicular to $Z_0$. For a fixed $X$, the latter vectors
can be chosen such that they are eigenvectors of the bilinear form
$\<J_Z(X),J_Z(X)\>$, restricted to the space $Z_0^\perp$. 
By using this
basis, the fourth term can be 
reduced to the third one in the first line of
(4.6).

The above formulas concern groups defined by general endomorphism 
spaces. In this paper we deal with groups defined by $ESW_A$'s and
we should explicitly compute 
the scalar curvature $\wt\kappa$ at points $(X,Z)$,
where $Z$ is an anticommutator. The endomorphisms defined by the
elements of the above introduced orthonormal basis
$e_0,e_1,\dots ,e_{l-1}$ on the Z-space ($e_0=Z/|Z|$) are 
denoted by $J_i$. Then the operators
$L_0=J_0^2$ and $L_\perp =\sum_{i=1}^{l-1}J_i^2$ commute and a common
eigensubspace decomposition can be established for them. In the 
following the scalar curvature 
$\wt\kappa$ is explicitly computed at a particular point $(X,Z)$
where $X$ is in the common eigensubspace of $L_0$ and $L_\perp$.
The corresponding eigenvalues are denoted by 
$\lambda_0$ and $\lambda_\perp$.
 
To be more precise, the scalar curvature will be computed at the points
of a 2-dimensional surface, called {\it Hopf hull}, which 
are included in a higher dimensional, so called, {\it X-hulls}.
These hulls are constructed as
follows.

For the unit vector $Z_0$, consider the one-parametric
family 
$S_X(sZ_0)$ 
of the X-crossections describing the so called
X-hull 
around $Z_0$. 
This X-hull is denoted by
$Hull_X(sZ_0)$. 
We can construct it by
an appropriate rotation of the graph of the function $|X|=\delta (|s|)$.
This construction shows that an X-hull is a k-dimensional 
manifold diffeomorphic to a sphere and also the s-parameter lines
on this surface are
well defined by the rotated graph. The point on the hull satisfying
$|X|=\delta (|S|)=0$ is the so called {\it vertex} of the hull. 
The sphere $S_X(0)$ in the middle is the {\it eye} of this 
manifold. This eye is shared by all of the X-hulls. It is a 
total-geodesic submanifold since it is fixed by the isometry
$(X,Z)\to (X,-Z)$. The
vertexes of the hulls form the so called {\it rim} of the sphere-type
manifold. This rim is nothing but the Z-sphere $S_Z(0)$ over
the origin of the X-space.

The Hopf hulls are sub-hulls of X-hulls, constructed as follows.
 
On the X-hull consider the
vector field $J_Z(X)$ tangent to the X-spheres. The integral curves
of this vector field are Euclidean circles (called {\it Hopf circles})
through an $X$ if and only if $X$ is an eigenvector of $J_Z^2$
with eigenvalue, say $\lambda_0$.
For other X-vectors these curves maybe not closed or they
are proper Euclidean ellipses.
If we fix a Hopf
circle $HC(0)$ on the sphere $S_X(0)$ at the origin and we consider
the s-parameter lines only through this circle, we get a 2-dimensional,
so called {\it Hopf hull},  
\,$H\!Hull_C(sZ_0)$, which is built up by a 1-parametric
family of parallel Hopf circles $HC(s)$. 

One can get this Hopf hull
by cutting it out from the ambient X-hull by the 3-dimensional
space $T_X$ spanned by the vectors $\{X\, ,\, J_Z(X)\, ,\, Z\}$.
From (1.3) and 
$$
[X,J_Z(X)]=\sum_\al\<J_Z(X),J_\al(X)\>Z_\al
\tag 4.7
$$
we get that $T_X$ is a total-geodesic manifold on the ambient space, 
for any $X$,
if and only if the $J_Z$ is an anticommutator. Then the $T_X$ is a
scaled metric Heisenberg group such that 
$|J_Z(X)|=\sqrt{-\lambda_0}|Z||X|$ holds. 

 Thus we get

\proclaim {Lemma 4.1} A Hopf hull, \, $H\!Hull_C(sZ_0)$,
is total geodesic
on a sphere-type manifold $\pa D$ for any Hopf-circle $C$ 
if and only if $Z_0$ is an 
anticommutator. These Hopf hulls are intersections of $\pa D$ by the
total-geodesic scaled Heisenberg groups $T_X$.
\endproclaim 

The scalar curvature $\wt\kappa$ (cf. above (4.5)) 
is computed on such a Hopf hull by formula
$$
\wt\kappa =\kappa -2Ricc(\bmu ,\bmu)+(TrB)^2-Tr(B^2)
\tag 4.8
$$
We use the new parameterization $\tau =|Z|^2=s^2$ on the
parameter lines. Then $D^\prime$ means differentiation 
with respect to this
variable.
By (1.16) and (4.1)-(4.5) one gets by a lengthy but
straightforward computation that the scalar curvature has the rational
form
$$
\wt\kappa ={Pol(\tau ,D(\tau ),D^\prime,D^{\prime\prime})\over
(4-\lambda_0t(D^\prime)^2)^2
(4D-\lambda_0\tau D(D^\prime)^2+4\tau )},
\tag 4.9
$$
where the coefficients of the polynomial $Pol$ (depending on 
$\lambda_0\, ,\,\lambda_\perp$ and on constants such as $k$ and $l$) 
can be computed by the following
formulas:
$$
\aligned
\kappa &={1\over4}(\lambda_0+Tr\,L_{\perp}),\\
-2Ricc&(\bmu ,\bmu )=C^2(-4D\lambda_\perp + \\ &\qquad
\qquad\lambda_0\tau
(-4D 
+(D^\prime)^2(2+D(\lambda_\perp +\lambda_0)))),\\
Tr\,B&=C(2(k-1-D^\prime(l-1))\\
&\qquad\qquad+D({1\over2}\lambda_\perp D^\prime +
\lambda_0({1\over2}D^\prime +\tau D^{\prime\prime})\Omega)), 
\\
-Tr\,B^2&=-C^2(4(k-1+(D^\prime )^2(l-1))
\\&\qquad +
\lambda_\perp D(2(D^\prime -(1+D^\prime )^2)+{1\over2}\lambda_0\tau )
\\
&\qquad
+4\lambda_0D({1\over2}D^\prime
+tD^{\prime\prime})\Omega
+(\lambda_0D({1\over2}D^\prime
+\tau D^{\prime\prime})\Omega)^2,
\\
\Omega&=4(4-\lambda_0\tau (D^\prime)^2)^{-1}\, ,\,
C^2 =(4D-\lambda_0\tau D(D^\prime)^2+4\tau )^{-1}.
\endaligned
\tag 4.10
$$

The explicit computation of the coefficients of the polynomial $Pol$
in formula (4.9) requires further tedious computations even in the
simple cases when the sphere-type
domain is nothing but the Euclidean sphere described by the equation
$|X|^2+|Z|^2=R^2$. In these cases $D(\tau )=R^2-\tau \, ,\,D^\prime =-1,
D^{\prime\prime}=0$
hold and $Pol$ is a fourth order polynomial of $t$. If $D(\tau )$ is a
higher order polynomial of $t$, then also $Pol(\tau )$ is a higher order
polynomial and, except only one polynomial, 
the $\wt\kappa (\tau )$ 
is a non-constant
rational function of $\tau$ and $\wt{\kappa}^\prime$ has only finite
many zero places. These examples show the wide range of the sphere
type domains for which $\wt{\kappa}^\prime\not =0$ almost everywhere.

These formulas allow to compare the inner scalar curvature 
$\wt{\kappa}_H$ of a Hopf hull with the scalar curvature $\wt\kappa$
of the ambient sphere type manifold. The above formulas can be applied
also for computing $\wt{\kappa}_H$ by the substitutions $L_\perp =0,
\lambda_\perp =0, k=2, l=1$. Then we get that the $\wt{\kappa}_H$ can
be expressed by means of the functions
$\tau ,D(\tau ),D^\prime,D^{\prime\prime}$ and by 
$\lambda_0$. The {\it Hopf curvature} $\wt{\kappa}_{HD}(\tau )$ of a
function $D(\tau)$ defining a sphere type domain is defined by the
scalar curvature $\wt{\kappa}_H$ such that $\lambda_0=-1$. I. e., it is
the scalar curvature of the sphere type domain defined by $D(\tau)$ on
the standard 3-dimensional Heisenberg group.

The scalar curvature has a simple form on groups with
Heisenberg-type endomorphism 
spaces.
(Reminder: An endomorphism space 
$J_{\bold z}$ is said to Heisenberg-type if
$J_Z^2=-|Z|^2id\, ,\,\forall Z\in
\bold z$.) In fact, by (4.6), the 
$\wt \kappa$ depends only on the functions
$\tau\, ,\, D\, ,\, D^\prime$ in these cases. From the above arguments 
we get that 
$\wt{\kappa}^\prime$ 
is non-vanishing almost everywhere
if and only if  
$\wt{\kappa}_{HD}^\prime$
is non-vanishing almost everywhere. 

\head\leftline{The Extension Theorem}
\endhead

The main result of this section is: 

\proclaim{Theorem 4.2} 
Let $\wt\Phi$ be an isometry between two sphere-type
hypersurfaces 
defined by the same function $D(|Z|^2)$ on the 
Heisenberg-type groups $N$ and $N^*$ 
such that the derivative 
$\wt{\kappa}_{HD}^\prime$ is
non-vanishing almost everywhere. (The abundance of such manifolds is
described below formula (4.10).)
Then the $\wt\Phi$ extends into an isometry of the form 
$\Phi =(\Phi_{/X},\Phi_{/Z})$ 
between the ambient
spaces, where the component maps are appropriate 
orthogonal transformations
on the X- resp. Z-spaces. 
Therefore the metrics $\wt g$ and $\wt{g}^*$ are isometric
if and only if the ambient spaces are isometric.
\endproclaim

We consider this problem first on groups defined by
Heisenberg-type $ESW_A$'s.
This simplification provides by good 
ideas for how to prove the theorem in the much more
complicated general cases.

\head\leftline{\it The Extension Theorem on Heisenberg-type groups} 
\endhead

In this case we consider a whole X-hull along with the
restriction of the scalar 
curvature $\wt\kappa$ of the sphere-type hypersurface onto it.
Since this scalar curvature 
depends only on the functions
$|Z|\, ,\, D\, ,\, D^\prime$, 
the vectors $grad\, \wt\kappa$ are always
pointing into the directions of the 
s-parameter lines. By the above remark, made about comparing 
$\wt{\kappa}^\prime$ and 
$\wt{\kappa}_{HD}^\prime$, we get 
that also  
$\wt{\kappa}^\prime$ is 
non-vanishing almost everywhere. Therefore $grad(\wt\kappa)$ 
is non-vanishing almost
everywhere on every X-hull. 
Since this vector field is invariant by isometries, 
the isometries must keep also
the parameter lines.
The vertex points are intersections of parameter-lines, therefore 
X-hulls are mapped to X-hulls such that vertex is mapped to vertex
and eye is mapped to eye. I. e. the isometries must keep
the X-crossections $S_X(Z)$ as well as the Z-crossections $S_Z(X)$.
 
From the construction of the X-hull it is clear that the s-parameter
lines identify the X-crossections 
$S_X(sZ_0)$ and, with respect to this identification,
an isometry $\wt\Phi$ defines the same map, 
$\wt\Phi_{/X}(Z_0)$, on the distinct X-spheres since all these maps
are identified with the one defined on $S_X(Z_0)$. 
In the following we show
that the $\wt\Phi_{/X}(Z_0)$ 
is the restriction of an orthogonal transformation (defined
on the ambient space) to the considered sphere. 

To prove this statement,
we explicitly compute 
the metric tensor
$$
g_{ij}=g\big(\pa_i,
\pa_j\big)\quad ,\quad g_{i\alpha}=g\big( \pa_i,
\pa_\alpha\big)\quad ,\quad g_{\alpha\beta}=g\big(\pa_\alpha
,\pa_\beta\big)
$$
on the ambient space. From (1.2) we get 
$$
\aligned
g_{ij}&=\delta_{ij} + \frac 1 {4} \<[X,\pa_i],[X, 
\pa_j]\>\\
&=\delta_{ij} + \frac 1 {4} D\sum_{\alpha =1}^l
\<\bold J_\alpha\big(X_0\big),\pa_i\>\<\bold J_\alpha
\big(X_0\big),\pa_j\>;
\\
g_{i\alpha}&=-\frac 1 {2} \<\bold J_\alpha\big(X\big),
\pa_i\> \quad ;\quad g_{\alpha\beta}=\delta_{\alpha \beta},
\endaligned
\tag 4.11
$$

The metric tensor $\delta_{ij}$ (in the above formula concerning
$g_{ij}$) defines the standard round metric on the above considered
X-spheres. Notice too that, because of the function 
$D$ in the second term,
the metrics on the X-spheres with different radius are non-homotetic.
Since the $\wt\Phi_{/X}(Z_0)$ keeps all of these 
different metrics, it must
keep them separately. I. e. it keeps
$\wt\delta_{ij}$ as well as $\sum
\<\bold J_\alpha\big(X_0\big),\wt\pa_i\>\<\bold J_\alpha
\big(X_0\big),\wt\pa_j\>$. Therefore it is derived from
an orthogonal transformation $\Phi_{/X}(Z_0)$ on the ambient X-space.
Actually, the transformations $\wt\Phi_{/X}(Z_0)$ and
$\Phi_{/X}(Z_0)$ do not depend on $Z_0$, since all the X-hulls share the
X-sphere $S_X(0)$ over $Z=0$ and isometries must keep this X-sphere at
the origin. Thus they can be simply 
denoted by $\wt{\Phi}_{/X}$ resp. $\Phi_{/X}$.

In the following we describe the isometry $\wt\Phi_Z(X)$ induced on
the Z-sphere $S_Z(X)$ over $X$. By the last formula of (4.11), each 
map $\wt{\Phi}_{/Z}(X):S_Z(X)\to S_{Z^*}(X^*)$
is the restriction of the orthogonal map
$\Phi_{/Z}(X):Z_X\to Z^*_{X^*}$, where $Z_X$ denotes the Z-space
over $X$.

On the other hand, from the formula given for $g_{i\alpha}$ in (4.11) we
get
$$
\<J_{Z^*}(X^* ),Y^*\>=\< J_Z(X),Y\> \quad ,\quad
J_Z=\Phi_{/X}^{-1}J_{Z^*}\Phi_{/X} ,
\tag 4.12
$$
i. e. the orthogonal transformations $\Phi_{/Z}(X)$ do not depend on $X$
and all are equal to the orthogonal transformation defined by the
conjugation
$J_{Z^*}=\Phi_{/X}J_{Z}\Phi_{/X}^{-1}$. 

This proves the desired 
extension theorem on Heisenberg-type groups completely.

\head\leftline{\it The extension in the general cases}
\endhead

Next we prove the extension theorem on spaces which are defined by
general $ESW_A$'s. The key idea of this general proof is the
same as before; first we show that any isometry, $\wt\Phi$, between
two sphere type domains defined by the same function $D(|Z|^2)$
leaves the X-space invariant such that it is the restriction of a 
uniquely determined orthogonal transformation $\Phi_{/X}:\bold v\to\bold
{v}^*$. This is the significant part of the proof  
since the construction of the appropriate orthogonal 
transformation $\Phi_{/Z}:\bold z\to\bold z^*$ 
can be completed on the same way than it is done above. 

The above technique used for constructing $\Phi_{/X}$ 
on Heisenberg type groups can not be directly applied
in the general cases because 
the hull should be
around an anticommutator and even on such an X-hull
the vector field $grad(\wt\kappa )$ is
not tangent to the s-parameter lines in general.
The latter tangent property is valid only 
on the much thiner 
Hopf hulls which are in a common eigensubspace of the commuting
operators $L_0$ and $L_\perp$. On these Hopf hulls both
$grad(\wt\kappa)$ and
$grad(\wt{\kappa}_H)$ are tangent to the s-parameter lines. (By the
argument explained at establishing the Hopf curvature 
$\wt{\kappa}_{HD}$
we get that these gradients are non-vanishing almost everywhere if and
only if 
$\wt{\kappa}_{HD}^\prime\not =0$ almost everywhere.) On the Hopf hulls
whose Hopf circles are still in an eigensubspace of $L_0$ but they are
not in an eigensubspace of $L_\perp$, only
$grad(\wt{\kappa}_H)$ 
is tangent to the parameter lines. If a parameter line is not inside of
an eigensubspace of $L_0$ then none of the above gradients is 
tangent to it.

We can deal with these difficulties in the following way.

Consider an X-hull around an anticommutator $Z$. By Lemma 4.1 and by
the technique developed for Heisenberg type groups we get that the
vertex $Z$ is mapped by the isometry $\wt\Phi$ to an anticommutator
$Z^*$ such that the Hopf hulls, which are inside
of a common eigensubspace
parametrized by the pair $(\lambda_0,\lambda_\perp )$ of eigenvalues,
are mapped to Hopf hulls which are inside of a common eigensubspace 
parametrized by the same pair of eigenvalues. On these Hopf hulls
also the parameter lines (together with parameterization) are kept by 
the isometry. Using all the total geodesic Hopf hulls whose Hopf circles
are in an eigensubspace $E_{\lambda_0}$ of $L_0$, the above statement
can be established for those parameter lines which are in 
$E_{\lambda_0}\oplus Z$. This means that the sub-eye 
$E_{\lambda_0}\cap S_X(0)$ is mapped to the sub-eye 
$E^*_{\lambda_0}\cap S^*_{X^*}(0)$.
It is clear too that the isometry restricted to
$(E_{\lambda_0}\cap S_X(0))\oplus Z$ 
extends to the sub-ambient space
$E_{\lambda_0}\oplus Z$ and thus it defines 
a uniquely determined orthogonal
transformation 
$\Phi_{\lambda_0}:\, E_{\lambda_0}\to E^*_{\lambda_0}$.
Our goal is to show that the orthogonal transformation 
$\Phi =\oplus\Phi_{\lambda_0}$, defined on the whole X-space, is the
desired one to the problem. To establish this statement it is enough
to prove that a general parameter line is mapped to parameter
line such that also the parameterization is kept on it.

The above defined orthogonal transformation $\Phi$ defines a
correspondence a\-mong the parameter lines of the two manifolds.
Next we show that the $\wt\Phi$ maps the corresponding parameter lines 
to each other.

A parameter line $p(s)$ can be represented in the form
$\delta X_0+sZ_0$. By (1.3), (4.1) and (4.3), the curvature
vector $\wt{\nabla}_PP$, where  
$P(\tau )=p^\prime (\tau)$ is the tangent vector, has the form
$$
\wt{\nabla}_PP =(\delta^{\prime\prime}/\delta^\prime )(P-Z_0)
-\delta^{\prime}J_{Z_0}(X_0).
$$  

This means that the parameter lines corresponding to each other 
are the solutions of the
same second order differential equation. This equation is invariant
under the action of isometries, therefore parameter lines are mapped
to parameter lines.
 
Since the complete eye is a total geodesic submanifold, it is mapped
to the eye 
$S^*_{X^*}(0)$. Furthermore, also the vertex-to-vertex
property is satisfied. Therefore
a general parameter line is mapped to parameter
line such that also the parameterization is kept on it.
 
It follows that isometries keep the X-spheres. More precisely,
if an X-sphere $S_X(\wt Z)$ is mapped to an $X^*$-sphere over
$\wt{Z}^*$ then the isometry can be described by the pair
$(\wt{\Phi}_{/X}(\wt X)=\wt{X}^*,
\wt{\Phi}_{/Z}(\wt Z)=\wt{Z}^* )$ such that both are the restrictions
of the corresponding orthogonal transformations from the pair
$(\Phi_{/X},\Phi_{/Z})$. By the same argument applied on Heisenberg type
groups we get that
$J_{Z^*}=\Phi_{/X}J_{Z}\Phi_{/X}^{-1}$, 
which proves the desired 
extension and the theorem in the general cases completely.

\head\leftline{\bf The Extension Theorem in the
solvable cases}\endhead

In the solvable cases we prove the non-isometry on such
sphere-type hypersurfaces which can be described by an equation of the
form $|X|^2=D(|Z|^2,t)=D(\tau ,t)$. By (3.24), the normal unit vector
of such a surface is
$$
\bmu^S =
Ct^{1\over2}(2X-D_\tau \bold J_Z(X))-
2CtD_\tau Z -cCtD_t\bold T=\bmu^S_X +\bmu^S_Z+\bmu^S_T,
\tag 4.13
$$
where $\bmu$ is considered to be 
an element of the Lie algebra,
furthermore
$$
C=\big (t(4|X|^2+ (D_\tau )^2|\bold J_Z(X)|^2)+t^2(
4\tau (D_\tau)^2+c^2D_t^2)\big)^{-\frac 1 {2}}.
\tag 4.14
$$

The unit vector $\bold t$ parallel to $\pa_t$
is perpendicular to $\bmu$ and
$\bmu_Z$, thus
$$
\bold t=(4D+(ctD_t)^2)^{-{1\over2}}(ctD_tX_0+2D^{1\over2}\bold T)=
\bold t_X+\bold t_T.
\tag 4.15
$$
Also this vector is considered to be a Lie algebra element.

These formulas together with (1.10), (1.11) and with
the formulas established on the nilpotent groups can be used for
computing the corresponding formulas also on the solvable groups.
Then we get:
$$
\aligned
M_S(\wt X_1,\wt X_2)&=t
M(\wt X_1,\wt X_2)+{1\over2}c^2CtD_t\<\wt X_1,\wt X_2\>\, ;\\
M_S(\wt X,\wt Z)&=
M_S(\wt Z,\wt X)=t^{1\over2}
M(\wt X,\wt Z)\, ;\\
M_S(\wt X,\bold t)&=
M_S(\bold t,\wt X)=P
\big (tD_t(1+{1\over4}c^2D_t)
\<X_0,\wt X\> \\
&\qquad -D\tau^{1\over2}D_{t\tau}\<J_{Z_0}(X_0),\wt X\>\big )\, ;
\\ P&=2
cCt(4D+(ctD_t)^2)^{-{1\over2}}\, ;\\
M_S(\wt Z_1,\wt Z_2)&=Ct(c^2-2D_\tau)\<\wt Z_1,\wt Z_2\>\, ;\\
M_S(\wt Z,\bold t)&=
M_S(\bold t,\wt Z)=P^*
\<J_{Z_0}(X_0),J_{\wt Z}(X_0)\>\, ;\\
P^*&=
{1\over2}cC
(4D+(ctD_t)^2)^{-{1\over2}}
t^{3\over2}\tau^{1\over2}D^{1\over2}D_\tau D_t\, ;\\
M_S(\bold t,\bold t)&=C\big (2t(1+{c^2\over4}D_t)|\bold t_X|^2 
-cDt^{1\over2}|\bold t_X||\bold t_Z| \\
&\qquad -c^2(tD_t+t^2D_{tt})|\bold t_T|^2\big ).
\endaligned
\tag 4.16
$$
(Also in this case
$D_\tau$ means differentiation with respect to
the argument $\tau =|Z|^2$.)

If $J_Z$ is an anticommutator,
the formulas for
$M_S(\wt X,\bold t)$ and
$M_S(\wt Z,\bold t)$
can be considerably simplified. In fact, the latter expression vanishes
and the first expression vanishes on the tangent vectors of the
form $\wt X=J_{\wt Z}(X)$ 
(these vectors are tangent to the surface,
since they are perpendicular to $\bmu^S_X$). Therefore it is non-trivial
only on the vectors $\wt X=J_Z(\bmu_X)$, furthermore, 
the plane spanned by this 
vector and by $\bold t$ is invariant by the action of the Weingarten
map $B_S$.

The proofs of the extension theorem 
can be straightforwardly adopted from the 
nilpotent
case to this solvable case.  

First choose
a unit vector $Z_0$ and consider the 
half-plane $\bold Z_0\oplus\bR_+$ parametrized by
$(s,t)$. Then the X-hull,\, $Hull_X(sZ_0,t)$,
is a 2-parametric family of
the X-spheres $S_X(sZ_0,t)$ defined on a closed domain diffeomorphic
to a closed disk. The boundary of this domain (which is diffeomorphic to
a circle) is the so called $(s,t)$-{\it rim}. The point $(s_v,1)$ resp.
$(0,t_v)$, where $D=|X|^2=0$, is called $Z$-{\it vertex}- resp. 
$t$-{\it vertex-parameters}. 
The corresponding points on the surface are the
corresponding vertexes. Also the Hopf hulls,\, $H\!Hull_C(sZ_0,t)$,
are 2-parametric families, $HC(s,t)$, of corresponding Hopf circles. 
One can get such a Hopf hull
by cutting it out from the ambient X-hull by the 4-dimensional
space $T_S$ spanned by the vectors $\{X\, ,\, J_Z(X)\, ,\, Z\, ,T\}$.
From (1.11) and (4.7)
we get that the space $T_S$ is total-geodesic in the ambient space 
if and only if the $J_Z$ is an anticommutator.

The scalar curvature $\wt \kappa_S$ depends on the functions listed in
(4.6) and the parameter $t$. 
Thus, by the very same arguments applied in the nilpotent case we 
get
\proclaim {Lemma 4.3} (A) 
A \, $H\!Hull_C(sZ_0,t)$ is total geodesic
on a Sphere-type manifold $\pa D$ if and only if $Z_0$ is an 
anticommutator. These Hopf hulls are the intersections of the sphere
type manifold by the total geodesic submanifolds $T_S$, which are the
solvable extensions of the corresponding Heisenberg subgroups, $T$,
introduced in the nilpotent case. I. e., the metrics on these
submanifolds are
complex hyperbolic metrics of constant holomorphic sectional 
curvature.
 
(B) If $Z_0$ is an anticommutator then 
$grad\,\wt\kappa_S(\not =0)$ is
tangent to the $(s,t)$-parameter plane
on a Hopf-Hull,  $H\!Hull_C(sZ_0,t)$, if and only if the Hopf circles
$HC(s,t)$ are in a common eigensubspace of the commuting operators
$L_0=J_0^2$ and $L_\perp =\sum_{i=1}^{l-1}J_i^2$.
\endproclaim 

The explicit computation of the scalar curvature can be performed
on a $H\!Hull_C$ by
using (4.8),(1.15),(1.16) and (4.10).
These lengthy computations are relatively 
simple when $D$ is a polynomial of $\tau$ and $t$.
(For instance, for Euclidean spheres with  center 
$(0,0,t_0)$ and radius $R$, this function is $D=R^2-\tau -(t-t_0)^2$.)
In these cases the $\wt \kappa_S$ is a rational function of $\tau$ and
$t$. 

The Hopf curvature $\wt{\kappa}_{HD}(s,t)$ is defined by the scalar
curvature of the sphere type manifold defined by $D$ on 
the standard complex hyperbolic space of $-1$ holomorphic sectional
curvature.
 
As in the nilpotent case we get

\proclaim{Theorem 4.4} Any isometry $\wt\Phi$ between two sphere type
manifolds, defined by the same function $D(s,t)$ such that 
$grad(\wt{\kappa}_{HD})\not =0$ almost everywhere, extends to
an isometry $\Phi$ between the ambient spaces $SN$ and $SN^*$.
I. e. $\wt g_c$ and $\wt g^*_c$
are isometric if and only if the ambient groups are isometric.
\endproclaim

\head\leftline
{\bf Extension- and nonisometry-theorems on the 
striking examples}
\endhead

The above proof of the extension theorem breaks
down on important hypersurfaces such as the geodesic spheres on the 
solvable groups
$SH^{(a,b)}_3$.
In \cite{Sz5} we
pointed out that the most striking examples can be constructed exactly
on these geodesic spheres.
In fact,
these geodesic spheres with the same radius 
are isospectral on spaces with the same
$a+b$, yet the spheres  
belonging to the 2-point homogeneous
space
$SH^{(a+b,0)}_3$
is homogeneous while the others are locally inhomogeneous.

Next we
establish the extension theorem along with other 
non-isometry theorems also on the geodesic spheres.
The key idea is an explicit computation of the
eigensubspaces of the Ricci curvature. The invariance of these 
eigensubspace-distributions guarantees that both the X-spaces and
Z-spaces are invariant under the actions of isometries and they extend
into an isometry between the ambient spaces.
Also in this section the nilpotent and the solvable cases are considered
separately

\head\leftline{ 
Extensions from the geodesic spheres of \,\,
$H^{(a,b)}_3$}
\endhead
 
We use the notations introduced in (2.12)-(2.14). In this case
$\{e_1=\bold i\, ,\,e_2=\bold j\, ,\, e_3=\bold k\}$
is a basis in the space $\bR^3$ of the imaginary quaternions and
$J_c=J_{e_c}:\bR^4\to\bR^4$ is defined by the appropriate left product
on the space $\bR^4$ of quaternionic numbers. The endomorphism
$J^{(a,b)}_c$ acting on $\bR^{4(a+b)}$ is introduced in (2.14).
The endomorphism spaces 
$J_{\bold z}^{(a,0)}$ and
$J_{\bold z}^{(0,b)}$
are used accordingly. Let us note that
$J_{\bold z}^{(a,0)}$ and
$J_{\bold z}^{(0,a)}$
are equivalent endomorphism spaces (in the sense of (1.7)) and they 
correspond to the left- resp. to the right-representation of $so(3)$
on $\bR^{4a}$.

We introduce also the distribution 
$\rho^{(a,b)}$
tangent to the X-spheres $S_X(Z)$
of the spaces $H^{(a,b)}_3$, spanned by the vectors 
$J_{\bold z}^{(a,b)}(X)$ at a vector $X$. Let us point out again that 
this distribution is considered as a regular X-distribution and the 
spanning vectors are regular X-vectors.
Therefore the $\rho$ is not perpendicular to the distribution $\wt z$
defined by the Z-vectors tangent to the considered surface at a point.
We introduce also the distribution $K^{(a,b)}$ consisting of vectors
perpendicular to $\rho^{(a,b)}\oplus \wt z$. 
From (1.2) we immediately get
that also this latter distribution is spanned by regular 
X-vectors, i. e.
$\rho^{(a,b)}\oplus K^{(a,b)}$ is an orthogonal direct sum decomposition
of
the tangent space of the Euclidean X-spheres around the origin.  

On the space $H^{(a,0)}_3$ (resp. on
$H^{(0,b)}_3$)
the distribution $\rho^{(a,0)}$ 
(resp. $\rho^{(0,b)}$) 
is integrable
and the 3-dimensional integral manifolds are the fibres of
a principal fibre
bundle with the structure group $SO(3)$. 
This fibration is nothing but
the quaternionic Hopf fibration and the factor space is the
2-point homogeneous quaternionic projective space \cite {Be}.
If $a\, >\,1$
(resp. $b\, >\,1$),
the distribution $K^{(a,0)}$ (resp. $K^{(0,b)}$)
is an irreducible connection
on this bundle with an irreducible curvature form $\omega (X,Y)=
[X,Y]_\rho ;X,Y\in K$.
This proves that 
$[K^{(a,b)},K^{(a,b)}]_\rho =\rho^{(a,b)}$. 
Thus we have 
\proclaim {Lemma 4.5} If $a,b\, >\, 1$, then
$[K^{(a,b)},K^{(a,b)}]_\rho =\rho^{(a,b)}$ 
and therefore the $K$ generates the whole tangent
space on a sphere $S_X(Z)$ by Lie brackets.
\endproclaim
 
In the following step we
compute the Ricci curvature on $\pa D$ and it turns out that both
$\rho\oplus \wt z$ 
and $K$ are eigensubspaces of this
Ricci operator with a completely different set of eigenvalues.
This observation offers more options for establishing the extension
theorem.

We use a special basis to compute the matrix of the Ricci curvature.
At a fixed point $(X,Z)$ on the 
hypersurface the unit normal vector $Z_0$
is denoted by $\bold i$, furthermore, the last two vectors from the
right handed orthonormal system $\{\bold i,\bold j,\bold k\}$ are
chosen such that they are tangents to the hypersurface. The unit vectors
$$
\gathered
\wt E_{\bold i}=J_{\bold i}(\bmu_{X0})=
-(4+|Z|^2(D^\prime )^2)^{-{1\over2}}(D^\prime |Z|X_0+2J_{\bold i}(X_0))
,\\  \wt E_{\bold j}=J_{\bold j}(X_0)\quad ,\quad \wt E_{\bold k}=
J_{\bold k}(X_0),
\endgathered
\tag 4.17
$$
(considered as Lie algebra elements)
are tangent to the hypersurface, standing perpendicular to the Z-space
and to the distribution 
$\rho^{(a,b)}$. 
Notice that 
$\wt E_{\bold j}$ and 
$\wt E_{\bold k}$ are tangent to
$\rho^{(a,b)}$ while the vector
$\wt E_{\bold i}$ is not tangent to this distribution,
expressing the fact
that $\rho$ and $\wt z$ are not perpendicular in general. We consider
an orthonormal basis $\{\wt K_1,\dots ,\wt K_{k-3}\}$ 
also on $K^{(a,b)}$ 
and the matrix of the Ricci operator $\wt r$ is computed 
with respect to the basis $\{\wt K_1,\dots ,\wt K_{k-3},\wt E_{\bold i},
\wt E_{\bold j},\wt E_{\bold k},\bold j,\bold k\}$. 
For the computations we
use the following formulas:
$$
\aligned
J_{\bold i}(\bmu_X)&=D^{1\over2}C(4+|Z|^2(D^\prime )^2)^{1\over2}
\wt E_{\bold i},\\
J_{\bold j}(\bmu_X)&=D^{1\over2}C(2
\wt E_{\bold j}+D^\prime |Z|\wt E_{\bold k}),\\
J_{\bold k}(\bmu_X)&=D^{1\over2}C(2
\wt E_{\bold k}-D^\prime |Z|\wt E_{\bold j}).
\endaligned
\tag 4.18
$$

Then by (4.1)-(4.5) we get that this matrix is of the form
   $$\wt r=\pmatrix
     \epsilon I_K &
     0 & 0 & 0\\
     0 &(\epsilon +E_{{ll}})I_{l}& 0 & 0\\ 
     0 & 0 &
    (\epsilon 
+E_{{LL}})
I_{L} & 
     E_{{L}\wt z} \\
     0 & 0 &
     E_{\wt z L} &
    (\epsilon 
+E_{\wt z\wt z})
I_{\wt z}
     \endpmatrix,\tag 4.19$$
where $I_K$, $I_{\wt z}$ and $I_{L}$ (resp. $I_{l}$)
are unit matrices on
the spaces $K$, $\wt z$ and on the space ${L}$
 spanned by the vectors
$\bold j$ and $\bold k$ (resp. on
 the 1-dimensional space 
 $l$ 
spanned by $\bold i$),
furthermore
$$
\aligned
\epsilon &=-{3\over2}+2C\, Tr(B)-4C^2 
=-{3\over2}+2C^2(Tr(b)-2)= \\ &=
-{3\over2}+2C^2\big (2(k-2)-DD^\prime - d_0D
\Omega \big ) ; \\
E_{{ll}}&=C^2\big (4+(3D-4)\Omega^{-1}\\
&\qquad -d_0D
(2(k-3)+D(d_0-D^\prime ))\Omega +2(d_0D)^2\Omega^2\big ); \\
E_{{LL}} &= C^2\big (4+(6D-4)\Omega^{-1}-
{1\over2}DD^\prime 
(2(k-3)
-{1\over2}DD^\prime ) \\
&\qquad\qquad +
{1\over2}D^2D^\prime d_0\Omega 
\big ); \\
E_{\wt z\wt z}&= {k\over4}+{3\over2}-2C^2\big ({1\over2}D\Omega^{-1} \\
&\qquad -
(1+D^\prime )(2(k-2)-D^\prime (D-2) )-D(D^\prime +1)d_0\Omega )\big ), 
\endaligned
\tag 4.20
$$
where $\Omega =4(4+|Z|^2(D^\prime ))^{-1}$ is introduced in (4.12)
and $d_0={1\over2}D^\prime +D^{\prime\prime}$ is introduced in (4.3).
The $2\times 2$ matrices
$E_{{L}\wt z} =E_{\wt z{L}}$ 
have the following form:
$$
\pmatrix
A & -B \\
B & A
\endpmatrix,
$$
where the functions $A$ and $B$ are
$$
\aligned
A&=C^2D^{1\over2}
\big ({1\over8}(3C^{-1} -D^{1\over2})\Omega^{-1}+6C^{-1} +
D(1+D^\prime )d_0\Omega\\
&\quad -(2D^{1\over2}D^\prime (D^\prime +2)+
(1+D^\prime)(2(k-1)-DD^\prime ))\big ) ;\\
B&=
C^2D^{1\over2}C^{-1}D^\prime |Z| \big
(k+2-{1\over2}DD^\prime-{1\over2}Dd_0\Omega\big ).
\endaligned
\tag 4.20'
$$

From these formulas and from the characteristic equation $det(\wt r
-\lambda I)=0$ we get that the subspaces $K$ and $\rho\oplus\wt z$
are eigensubspaces of the Ricci operator and the  
eigenvalue $\epsilon$ on $K$ is different from the other eigenvalues
if and only if the following determinant
   $$ det\pmatrix
     E_{{LL}}I_{L} &
     E_{L\wt z} \\
     E_{\wt z L} &
     E_{\wt z\wt z} I_{\wt z}
     \endpmatrix
=(A^2+B^2-E_{LL}E_{\wt z\wt z})^2\tag 4.21$$
is non-zero. Since this determinant is zero on an open set only
in the case when
the function $D$ satisfies a certain differential equation on some
open intervals, it is clear that 
this last assumption is satisfied on an everywhere
dense open set for the most general sphere-type manifolds in 
$H_3^{(a,b)}$. Then on this set
the Ricci tensor has distinct eigenvalues on the
invariant subspaces $K$ and $\rho\oplus\wt z$. 
Since we concentrate on the geodesic spheres in the next section, 
we would like to give more details about certain particular cases
in order to prepare the next section.

Later we will see that the geodesic spheres on the solvable extensions
intersect the nilpotent level sets in Sphere-type
surfaces described by functions
of the form $D(\tau)=\sqrt{Q-\tau}+Q^*$, where $\tau =|Z|^2$ and $Q\, ,
Q^*$ are constants. In this case we introduce the new variable
$u=\sqrt{Q-\tau}$. Then all the functions $D\, ,\, D^\prime \, ,\, d_0\,
,\, C^2\, ,\,\Omega$ are rational functions of $u$, while the functions
$D^{1\over2}=\sqrt{u+Q^*}$ and
$ 
C^{-1}(u)
={1\over2u}G(u)
={1\over2u}H^{1\over2}(u),
$
where
$$
H(u)
=-16u^4+15u^3+(16Q+15Q^*)u^2-Qu-QQ^*
\tag 4.22
$$
are non-rational functions of $u$. (One can prove, by an 
elementary computation, that the polynomial $H(u)$ can be written
in the quadratic form $H(u)=(q_1u^2+q_2u+q_3)^2$ only for particular
constants $Q$ and $Q^*$, and even in this particular case the 
coefficients $q_s$ are pure imaginary numbers ($q_1=\pm{4\bold i}$
can be checked immediately)). This argument proves that the function
$B^2-E_{{LL}}E_{\wt z\wt z}$ is a rational function of $u$, however,
the function $A^2$ has the following non-rational form
$$
A^2(u)=R_1(u)+R_2(u)D^{1\over2}(u)+R_3(u)H^{1\over2}(u),
\tag 4.23
$$
where the functions $R_i(u)$ are non-trivial rational functions. This
proves that the determinant in (4.17) is non-vanishing almost 
everywhere since the non-rational terms can not be canceled out
from this function either. 

The proof of the following 
extension and non-isometry theorem is completely prepared
by the above considerations.

\proclaim{Theorem 4.6} Suppose that the function defined in (4.21) is
non-zero on an everywhere dense open set (this assumption is most
widely satisfied, including the manifolds
described above around formula (4.22)).
Then the isometries $\wt\Phi :\pa D\to\pa D^*$
keep both the X-spaces and the Z-spaces and they extend to
isometries $\Phi$ acting between the ambient spaces. Therefore
the metrics $\wt g$ and $\wt{g}^*$ are isometric if and only if
the ambient metrics $g$ and $g^*$ are isometric.
\endproclaim

\demo{Proof}
First we suppose that $dim(K)>0$. Since $[K,K]_\rho =\rho$ and the
$\wt \Phi$ keeps the distribution $K$ by the above arguments,
the $\wt \Phi$ keeps the tangent spaces of the X-spheres $S_X$.
Therefore the image of $S_X(Z)$ is an $X^*$-sphere $S_{X^*}
(Z^*)$. Thus the $\wt \Phi$ defines an orthogonal transformation
between these 2 spheres, transforming the vector field $J_Z(X)$ to $
J_{Z^*}(X^* )$.
(The last statement follows from (4.11) by the
same arguments used there, since also in this case 
the $\wt \Phi$ keeps the tensors $\delta_{ij}$ and
$\sum_\al \<J_\al (X),\wt\pa_i\>\<J_\al (X),\wt\pa_j\>$
together with the tensor $g_{ij}=-1/2\<J_j(X),\wt\pa_i\>$ 
separately. Therefore it is derived from an orthogonal transformation
on the ambient space such that 
also the form $\<J_0(X),\wt\pa_i\>$ is preserved. This latter statement
can be also proven by using the invariant property of the eigensubspace
$l$ of the Ricci operator $\wt r$.)

Let $\wt\rho$ be the subspace in $\rho\oplus \wt t$ perpendicular to
$\wt z$. Then the $\wt \Phi$ keeps $\wt\rho$ by (4.13) and
by the previous argument. Consequently,
it keeps also $\wt z$. I. e., a Z-sphere $S_Z(X)$ is mapped to the
$Z^*$-sphere $S_{Z^*}(X^{*})$ 
and the $\wt\Phi$
extends to an orthogonal transformation between the ambient spaces
$z_X$ and $z_{X^*}$. Since the ambient spaces are isometric 
if and only if $(a,b)=(a^* ,b^* )$ up to an order, this
extension theorem proves also the
non-isometry completely.
\enddemo

\noindent{\bf Remark 4.7}
The non-isometry can be established also on other ways as follows.

($\bold A$)
For the vector fields $U$ and $V$ tangent to
$\rho^{(a,b)}\oplus t^2$, 
let $L(U,V)$ be the orthogonal
projection of $[U,V]$ onto $K^{(a,b)}$. Then $L$ is obviously a tensor
field of type (2,1) on $\pa D$ such that
it is invariant with respect to the isometries
of the space. It turns out too that the $L$ vanishes exactly at the 
points of the form $(X^{(a)},Z)$ or $(X^{(b)},Z)$. 
I. e., the induced metrics
on the hypersurface $\pa D$ of the spaces
$H^{(a,b)}_3$ and
$H^{(a^*,\, b^*)}_3$ 
can not be isometric unless $(a,b)=(a^*,\,b^*)$ up to an order.

($\bold B$)
One can demonstrate the non-isometry also by 
determining the isometries on the considered
hypersurfaces. It turns out that the group of isometries
is 
$\{O(\bold H^a)\times O(\bold H^b)\} SO(3)$, 
where $O(\bold H^{c})$ is the 
quaternionic orthogonal group acting on $\bold H^{c}$. This also
proves the above local non-isometry on 
$\pa D$. This proof is not independent from the above extension
theorem, since the above isometry group is determined by
knowing that the isometries on $\pa D$ are the restrictions of 
those isometries on the ambient space which fix the origin.

($\bold C$) One can establish a general(!) non-isometry proof (not just
on $H^{(a,b)}_3$) also by \cite{Sz5,Proposition 5.4}, where the
isotonal property of the corresponding curvature oprators is established
(cf. also Proposition 1.3 in this paper). The complete 
proof on sphere-type
manifolds needs the extension theorem also in this case. 

\head\leftline{The non-isometry proofs on the geodesic
spheres of $SH^{(a,b)}_3$}
\endhead

In order to establish these results on the solvable extension $SN$,
first we explicitly 
compute the equation of a geodesic sphere around the origin
$(0,0,1)$. This computation can be carried out by using
the generalized Cayley transform constructed on the solvable extensions
of Heisenberg type groups \cite{CDKR}. This transform maps the
unit ball
$$
B=\{(X,Z,t)\in\bold n\oplus\bold a\, |\, |X|^2+|Z|^2+t^2=r^2<\,1\}
\tag 4.24
$$
onto $SN$, by the formula
$$
C(X,Z,t)=((1-t)^2 +|Z|^2 )^{-1} \big( 2(1-t+J_Z)(X)\, ,\, 2Z\, ,
\, 1-r^2 \big ).
\tag 4.25
$$
By pulling back we get the ball-representation of the metric,
having the property that the geodesics through 
the origin are nothing but
the rays
$(tanh(s)/r)(X,Z,t)$ such that they are
parameterized by the arc-length $s$ starting at the origin. I. e.,
the considered geodesic spheres with radius $s$ match with the Euclidean
spheres with radius $tanh(s)$ on the Ball-model. Then, by computing the 
inverse Cayley map (the reader can consult for more details with
\cite{CDKR} (pages 14-15)), by a routine computation we get
\proclaim{Lemma 4.8}
The equation of a geodesic sphere of radius $s$ around the origin
$(0,0,1)$ of $SN$ is
$$
|X|^2=4((e^s+e^{-s}+2)t-|Z|^2)^{1\over2}-4(t+1).
\tag 4.26
$$
\endproclaim

Notice that on groups
$SH^{(a,b)}_3$
with the same $a+b$, the geodesic spheres with the same radius 
$R$ around the
origin are the same level sets, described by
the same equation.

In the following step we compute the Ricci curvature $\wt r_S$ on the
geodesic spheres 
with respect to the basis $\{\wt K_1,\dots ,\wt K_{k-3},\wt E_{\bold i},
\wt E_{\bold j},\wt E_{\bold k},\bold j,\bold k,\bold t\}$. By
(4.16), (1.13), (4.1), (4.5), (4.13) and by
$$
\aligned
R(\bmu^S,\wt U,\wt V,\bmu^S)&=
R(\bmu_X^S+\bmu^S_Z,\wt U,\wt V,\bmu_X^S+\bmu_Z^S)+
R(\bmu_T^S,\wt U,\wt V,\bmu_T^S)\\ &\quad +
R(\bmu_T^S,
\wt U,\wt V,
\bmu_X^S+\bmu_Z^S)
+R(
\bmu_X^S+\bmu_Z^S,
\wt U,\wt V,
\bmu_T^S)
\endaligned
$$
 we get 
   $$\wt r_S=\pmatrix
     \sigma I_K &
     0 & 0 & 0 & 0\\
     0 &(\sigma +S_{{ll}})& 0 & 0 & S_{l\bold t}
\\ 
     0 & 0 &
    (\sigma 
+S_{{LL}})
I_{L} & t^{1\over2} 
     S_{{L}\wt z} &0 \\
     0 & 0 &t^{1\over2}
     S_{\wt z L} &
    (\sigma 
+S_{\wt z\wt z})
I_{\wt z} &0 \\
0&S_{\bold t l}&0&0& \sigma+S_{\bold t\bold t}
     \endpmatrix,\tag 4.27$$

By the very same proof given in the nilpotent case we get that the
eigenvalue $\sigma$ on the eigenspace $K$ is different from the other
eigenvalues on an everywhere dense open set and therefore this 
distribution is invariant by the actions of the isometries. 

In fact,
the determinant corresponding to (4.21) is non-vanishing also in this
case, since the corresponding term $A^2$ has now the form
$$
A^2(u,t)=R_1(u,t)+R_2(u,t)D^{1\over2}(u,t)+R_3(u,t)H^{1\over2}(u,t),
\tag 4.28
$$
where
$u=((e^R+e^{-R}+2)t-\tau )^{1\over2}$, 
while the other terms are rational
functions. A straitforward computation shows that also 
$S_{{ll}}S_{\bold{tt}}-S_{l\bold t}^2$ is non-vanishing 
on an everywhere dense open set, which proves the above statement
concerning the distinctness of $\sigma$ from the other 
eigenvalues.

The extension- and nonisometry-proofs can be similarly 
established than in the nilpotent case.

Alike to the nilpotent case, also the
invariant tensor field
$L(\wt U,\wt V)$ defined on the distribution $\rho\oplus\wt z\oplus
\bold t$ can be used for establishing the non-isomery.
This tensor field vanishes exactly at the X-vectors of
the form $(X^{(a)},0)$ or $(0,X^{(b)})$, proving the
non-isometry of the considered manifolds.

If $ab\not =0$,
the group of isometries is the 
non-transitive
group
$\{O(\bold H^a)\times O(\bold H^b)\} SO(3)$ on the geodesic spheres,
while the geodesic spheres of the 2-point homogeneous spaces 
$SH^{(a+b,0)}_3$ 
have transitive groups of isometries whose unit component 
is isomorphic to 
$\bold{Sp}(a+b+1)\bold{Sp}(1)$. 
This is the third proof of the non-isometry. This demonstration is
not independent from the extension theorem, since the 
isometry group is determined by the fact that the isometries on the
geodesic sphere are restrictions of those isometries on the ambient
space which fix the center of the sphere.

A different type of the non-isometry proofs can be 
established by Propsition 1.3, though also this proof involves
the extension theorem.
  
By summing up, both the isospectrality and the non-isometry theorems, 
we have
\proclaim {Cornucopia Theorem 4.9}
(A)Let $J_{\bold z}$ be an endomorphism space with anticommutator
(i. e., an $ESW_A$) such 
that it either contains a nonabelian Lie subalgebra
or it is one of the irreducible Cliffordian endomorphism spaces
$J_{4k+3}^{(1,0)}\simeq J_{4k+3}^{(0,1)}$. 
Consider the family, determined by the constant $(a+b)$, 
of 2-step 
nilpotent metric Lie groups 
$(N^{(a,b)}_J,g^{(a,b)})$
(resp. the family 
$(SN^{(a,b)}_J,g^{(a,b)})$
of solvable extensions) defined by the endomorphism spaces 
$J^{(a,b)}_{\bold z}$ (cf. (2.12)-(2.14)).
Such a family is represented on the same manifold
$M=\bold{R}^{k(a+b)+l}$
(resp. on $M=\bold{R}^{k(a+b)+l}\times\bold{R}_+$ in the solvable case).

While the induced metrics $\wt g^{(a,b)}$ and $\wt g^{(b,a)}$ are 
isometric, the other metrics from the family have different local
geometries than $\wt g^{(a,b)}$, 
on any sphere-type hypersurface $\pa D\subset M$ defined
by the same function $\varphi (|X|,|Z|)=0$
(resp. $\varphi (|X|,|Z|,t)=0$) where 
the condition $grad(\wt{\kappa}_{HD})\not =0$ is satisfied
almost everywhere on the corresponding Hopf-hull. Yet the metrics
$\wt g^{(a,b)}$ on a $\pa D$, 
belonging to a family, are isospectral.

(B)The above statement is established also on the
geodesic spheres
of the solvable extensions
$SH^{(a,b)}_3$ (the technique applied for proving 
Theorem (A) breaks down in this case).
I. e., the geodesic spheres having
the same radius and belonging to the same family have 
different local geometries unless $(a,b)=(a^*,\, b^*)$ up to an order.
Yet, the induced metrics are isospectral also in this case.

The geodesic spheres on  
$SH^{(a+b,0)}_3$
are homogeneous, while the geodesic spheres 
on the other manifolds 
$SH^{(a,b)}_3$
are locally
inhomogeneous. This demonstrates the fact:

"One can not hear the local homogeneity
property even on the most simple closed manifolds, namely, which 
are diffeomorphic to Euclidean
spheres." 
\endproclaim

The abundance of the isospectrality examples constructed in these papers
is due to the abundance of the $ESW_A$'s, described in Section 2, and to
the great variety of the sphere-type manifolds
which can be chosen for a fixed
family of endomorphism spaces, both on the nilpotent groups and the
solvable extensions. Let us mention again that the isospectrality
theorem is established for any sphere-type manifold and the non-isometry
theorems are established for the particular manifolds defined
by equations of the form 
$\varphi (|X|,|Z|)=0$
(resp. $\varphi (|X|,|Z|,t)=0$)
only because we have not wanted to make the proofs even more
complicated than they are in this simplified situation.
It is highly possible that the extension and non-isometry proofs can
be extended to sphere-type manifolds defined by general functions
of the form
$\varphi (|X|,Z)$
(resp. $\varphi (|X|,Z,t)$).

\enddocument